\documentclass[a4paper, reqno, 12pt]{amsart}

\usepackage{geometry}       

\usepackage{float}
\usepackage{bm}
\usepackage{fullpage,xcolor}
\usepackage[mathscr]{euscript}
\usepackage[all]{xy}
\usepackage{epsfig}
\usepackage{amsfonts}
\usepackage{mathptmx}

\usepackage{graphicx}
\usepackage{amssymb} 
\usepackage{amsmath}
\usepackage{amsthm}
\usepackage{mathrsfs}
\usepackage{epstopdf}
\usepackage{url}
\usepackage[msc-links,numeric]{amsrefs}
\usepackage{tikz}

\usepackage{color}
\makeatletter

\@addtoreset{equation}{section}
\makeatother 

\theoremstyle{plain}
\newtheorem{theorem}{Theorem}[section]

\newtheorem{proposition}[theorem]{Proposition}

\theoremstyle{definition}
\newtheorem{definition}[theorem]{Definition}
\newtheorem{example}[theorem]{Example}
\newtheorem{remark}[theorem]{Remark}

\usepackage{latexsym}

\begin{document}

\title[Non-Classical Knots]{Algebraic and Geometric Aspects of Non-Classical Knots}

\author{Ioannis Diamantis}
\address{Department of Data Analytics and Digitalisation,
School of Business and Economics, Maastricht University,
P.O.Box 616, 6200 MD, Maastricht,
The Netherlands.
}
\email{i.diamantis@maastrichtuniversity.nl}



\setcounter{section}{-1}

\date{}



\begin{abstract}
Non-classical knot theory refers to a family of extensions of classical knot theory in which one or more of the basic ingredients of the classical framework are modified. These modifications may affect the local structure of crossings, the allowed Reidemeister-type moves, the ambient space, the global form of diagrams, or the additional data carried by diagrams. Examples include skein modules of three-manifolds, pseudo knots, singular knots, stuck knots, bonded knots, tied links, virtual knots, welded knots, knotoids, and related braid-type structures. In this survey, we present an overview of several non-classical knot theories from a comparative geometric and algebraic perspective. We examine how changes in ambient topology, new crossing types, rigidity constraints, auxiliary relational data, modified isotopy relations, and open or virtual diagrammatic settings lead to generalized knot theories with distinct topological and combinatorial features. Particular emphasis is placed on braid-theoretic formulations, skein-theoretic methods, trace constructions, and extensions of classical polynomial invariants. We conclude by outlining open problems concerning generalized algebraic structures, hybrid diagrammatic theories, and the relationships among different non-classical knot theories.
\end{abstract}

\maketitle

\section{Introduction}

Classical knot theory studies embeddings of one or more disjoint circles in 3-dimensional space, considered up to ambient isotopy. Its diagrammatic formulation, originating in the work of Reidemeister~\cite{Reidemeister1932} and Alexander~\cite{Alexander1923}, translates geometric questions about embedded curves into combinatorial questions about planar diagrams and local moves. In this setting, two knot or link diagrams represent the same isotopy class if and only if they are related by a finite sequence of Reidemeister moves. This correspondence lies at the heart of classical knot theory and leads naturally to algebraic and combinatorial structures such as braid groups \cite{Artin1925,Birman1974}, skein relations \cite{Conway1970,Kauffman1987,PrzytyckiTraczyk1987}, quandles \cite{Joyce1982}, state-sum models \cite{Kauffman1987}, and polynomial invariants \cite{Alexander1928,Jones1985,FreydYetterHosteLickorishMillettOcneanu1985,
PrzytyckiTraczyk1987,Kauffman1987}.

The relationship between knots and braids plays a particularly important role. Alexander's theorem~\cite{Alexander1923} states that every oriented link may be represented as the closure of a braid, while Markov's theorem \cite{Markov1935,Birman1974} characterizes when two braids have isotopic closures. These results place braid groups at the algebraic core of knot theory and connect diagrammatic topology with representation theory and diagrammatic algebras such as the Hecke algebra, the Temperley--Lieb algebra, and the Birman--Murakami--Wenzl algebra \cite{Jones1987,BirmanWenzl1989,Murakami1987}. Through suitable representations and trace constructions, one obtains many of the fundamental polynomial invariants of classical knot theory, including the Jones, HOMFLYPT, and Kauffman polynomials
\cite{Jones1985,FreydYetterHosteLickorishMillettOcneanu1985,
PrzytyckiTraczyk1987,Kauffman1987,Murakami1987}.

At the same time, classical knot theory is based on a highly specific diagrammatic framework. Crossings are fully determined as overcrossings or undercrossings, the ambient space is usually \(S^3\) or \(\mathbb R^3\), the objects are closed curves, and the equivalence relation is generated by the classical Reidemeister moves. Many natural situations require extensions of this framework. Crossings may be ambiguous, singular, or constrained; diagrams may carry additional relational data; the ambient surface or three-manifold may change; or the underlying diagrammatic object may no longer be closed. These possibilities have led to a broad family of non-classical knot theories \cite{Kauffman1991KnotsPhysics}.

One major direction is obtained by changing the ambient three-manifold. Skein modules, introduced by Przytycki and Turaev~\cite{Przytycki1991,Turaev1990SkeinSolidTorus},
extend the idea of polynomial link invariants from links in \(S^3\) to links in arbitrary three-manifolds. In \(S^3\), a skein relation often determines a polynomial invariant. In a general three-manifold, the same local relation usually produces a module generated by links in the manifold. This shift from polynomials to modules provides a powerful framework for studying knots and links in spaces such as the solid torus, lens spaces, handlebodies, and \(S^1\times S^2\) \cite{Turaev1990SkeinSolidTorus,HaringOldenburgLambropoulou2002Handlebodies,
Diamantis2019AlternativeBasisST,Diamantis2019KBSMGenus2Handlebody,
Diamantis2024KBSMS1S2,Diamantis2025HOMFLYPTS1S2}, and it connects naturally with braid-theoretic methods through mixed braid groups and generalized Hecke-type algebras \cite{Lambropoulou1999HeckeTypeB,LambropoulouRourke1997Markov3Manifolds,
DiamantisLambropoulou2015BraidEquivalence3Manifolds,
Diamantis2025SurveySkeinModules}.

Another major direction consists of modifying the local structure of crossings. Singular knot theory introduces prescribed transverse double points and plays a central role in the theory of Vassiliev finite-type invariants \cite{Vassiliev1990,Baez1992,BarNatan1995}. Pseudo knot theory replaces some classical crossings by unresolved crossings, leading to diagrammatic models with incomplete crossing information
\cite{Hanaki2010,Dye2010,HenrichEtAl2013Pseudoknots}. Stuck knot theory \cite{Diamantis2026StuckKnots} provides a different crossing-based generalization: the crossing remains classical, with specified over/under information, but its behavior under isotopy is constrained. In each of these theories, new local crossing data require corresponding modifications of Reidemeister moves, braid structures, skein relations, and polynomial invariants.

A further direction enriches diagrams with additional structure rather than modifying crossings themselves. Tied links introduce ties encoding partition-like relations among strands or components \cite{AicardiJuyumaya2016TiedLinks}, while bonded knots incorporate embedded auxiliary connections between points of a diagram \cite{GoundaroulisEtAl2017BondedKnotoids,DiamantisKauffmanLambropoulou2025Bonded}, motivated in part by applications in molecular and network-like systems  \cite{DorierEtAl2018KnotoidProtein,GoundaroulisEtAl2017BondedKnotoids}. These theories show that knot diagrams may be generalized not only through new crossing types or new ambient spaces, but also through additional geometric or combinatorial data attached to the diagram.

Other generalizations change the ambient or global diagrammatic setting. Virtual knot theory, introduced by Kauffman~\cite{Kauffman1999}, enlarges classical knot diagrams by allowing virtual crossings, which encode artifacts of representing knots in thickened surfaces by planar diagrams. This viewpoint connects knot theory with surface topology and Gauss diagrams
\cite{CarterKamadaSaito2002,Kuperberg2003}. Welded knots arise as a quotient of virtual knot theory obtained by allowing one of the forbidden moves, and are closely related to loop braid groups and ribbon structures in four-dimensional topology \cite{FennRimanyiRourke1997BraidPermutation}. Knotoids, introduced by Turaev~\cite{Turaev2012}, modify the global nature of the object itself by replacing closed curves with open-ended diagrams having two endpoints. This apparently small change produces a rich theory with new closure operations, equivalence relations, and applications, particularly in the study of open molecular chains
\cite{GugumcuLambropoulou2021Braidoids,
GoundaroulisEtAl2017BondedKnotoids,DorierEtAl2018KnotoidProtein}.

Many of these non-classical theories admit braid-theoretic formulations that extend the classical braid group framework. In several cases, braid groups are replaced by braid monoids \cite{BardakovJablanWang2016PseudoBraids,FennKeymanRourke1998SingularBraid}, mixed braid groups \cite{LambropoulouRourke1997Markov3Manifolds,
DiamantisLambropoulou2015BraidEquivalence3Manifolds}, virtual or welded braid groups \cite{KauffmanLambropoulou2006VirtualBraids,FennRimanyiRourke1997BraidPermutation}, braidoids \cite{GugumcuLambropoulou2021Braidoids}, or enriched braid structures containing additional generators corresponding to generalized crossings, bonds, ties, or other diagrammatic data \cite{BardakovJablanWang2016PseudoBraids,
DiamantisKauffmanLambropoulou2025Bonded,
AicardiJuyumaya2016TiedLinks}. This leads naturally to generalized algebraic structures, including Hecke-type algebras, braid-and-tie algebras \cite{AicardiJuyumaya2000BraidsTies,
AicardiJuyumaya2016TiedLinks}, skein-theoretic quotients, state-sum models, and trace constructions adapted to non-classical settings \cite{Jones1987,ParisRabenda2004,Diamantis2026HOMFLYPTPseudoResolution}.

The purpose of the present paper is to provide a comparative survey of several major non-classical knot theories from geometric, diagrammatic, and algebraic perspectives. Rather than giving a complete technical account of every theory, we focus on the structural mechanisms through which the classical framework is extended. One begins with a generalized class of diagrams, modifies the corresponding local moves and equivalence relations, develops a braid-like formulation when possible, and studies the associated algebraic structures and invariants. While this paradigm originates in classical knot theory, it becomes substantially richer in the non-classical setting, where ambiguity, degeneration, rigidity, bonding, virtuality, and ambient topology may all become part of the structure.

The paper is organized as follows. Section~2 recalls the classical framework that serves as a structural template: diagrams, Reidemeister moves, polynomial invariants, braid groups, and the braid--algebra--trace mechanism. Section~3 discusses skein modules and knots in three-manifolds, with emphasis on the solid torus, lens spaces, and braid-theoretic approaches. Section~4 treats pseudo and singular knot theories, emphasizing their parallel braid-theoretic and algebraic structures as well as their different geometric interpretations. Section~5 surveys stuck knots, where classical crossings are equipped with rigidity constraints and polynomial invariants are adapted to the restricted move system. Section~6 discusses bonded knots and tied links as examples of diagrammatic enrichments by auxiliary relational data. Section~7 collects further diagrammatic and ambient generalizations, including virtual knots, welded knots, knotoids, braidoids, and related variants. The paper concludes with open problems and directions for future work.


\section{Classical Knot Theory as a Structural Template}

Classical knot theory provides the geometric, diagrammatic, and algebraic foundation for many of the generalized theories discussed in this survey \cite{Reidemeister1932,Birman1974,Lickorish1997,Kauffman1991KnotsPhysics}. Its central insight is that topological information about embedded curves in three-dimensional space can be encoded by planar diagrams, provided that one keeps track of local crossing data and works modulo appropriate local moves. In this way, a geometric problem about embeddings becomes a combinatorial problem about diagrams \cite{Reidemeister1932,Lickorish1997}. The non-classical theories considered later arise by altering one or more ingredients of this classical framework: the nature of crossings, the allowed moves, the ambient space, or the additional structure carried by the diagram.

\subsection{Knots, links, and diagrams}

A \emph{knot} is a smooth embedding \(K:S^1\hookrightarrow S^3\). More generally, a \emph{link} is a smooth embedding of a finite disjoint union of circles, \(L:S^1\sqcup\cdots\sqcup S^1\hookrightarrow S^3\). The embedded circles are called the \emph{components} of the link, so that a knot is precisely a link with one component.

\begin{figure}[ht]
    \centering
    \includegraphics[width=0.25\textwidth]{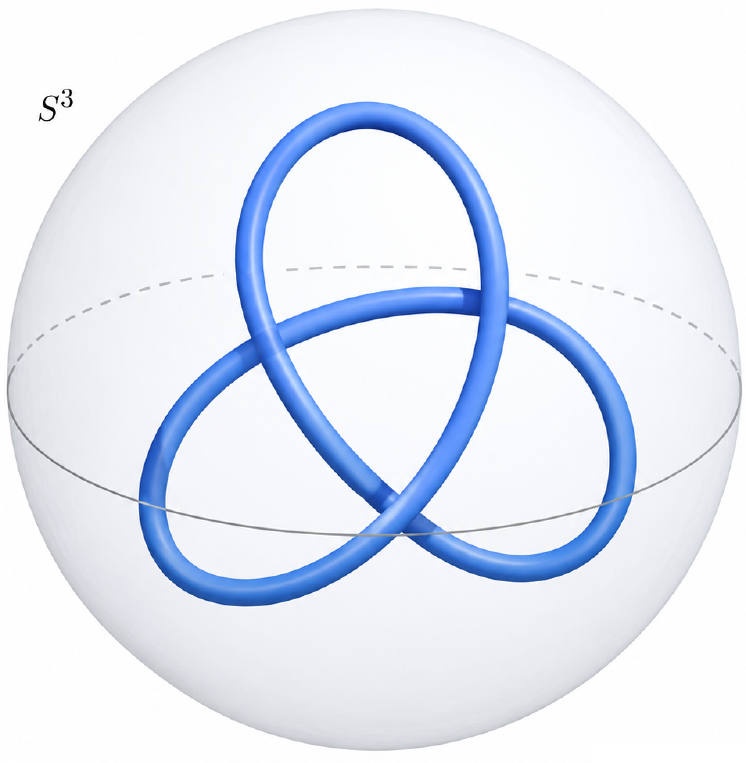}
    \caption{A schematic visualization of a knot in \(S^3\), drawn inside a three-ball after identifying \(S^3\) with the one-point compactification of \(\mathbb R^3\).}
    \label{fig:classical-knot-S3}
\end{figure}

Since \(S^3\) may be viewed as the one-point compactification of \(\mathbb R^3\), one often represents a knot by its image in \(\mathbb R^3\), with the point at infinity understood implicitly. In figures, this image is often drawn inside a three-ball, as in Figure~\ref{fig:classical-knot-S3}.

Knots are regarded as flexible embedded curves rather than rigid geometric objects. Thus, their precise shape in space is not important; what matters is
whether one embedded circle can be continuously deformed into another without cutting the curve or allowing it to pass through itself. This idea is modeled mathematically by \emph{ambient isotopy}. More precisely, two embeddings \(K_0,K_1:S^1\hookrightarrow S^3\) are considered equivalent if there exists a continuous family of homeomorphisms
\[
H_t:S^3\to S^3,\qquad t\in[0,1],
\]
with \(H_0\) equal to the identity and \(H_1(K_0)=K_1\). This formalizes the idea that one knot may be deformed into another by moving the surrounding space, without cutting the curve or allowing it to pass through itself.

\begin{figure}[ht]
    \centering
    \includegraphics[width=0.65\textwidth]{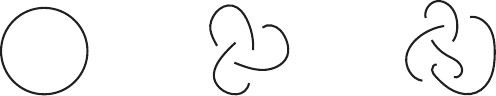}
    \caption{Basic examples of classical knots: the unknot, the trefoil knot, and the figure-eight knot.}
    \label{fig:basic-classical-knots}
\end{figure}

Although knots and links are defined as embedded objects in three-dimensional space, they are usually studied through {\it planar diagrams}. Such diagrams are obtained by projecting the embedded curves generically onto a plane. We note that a generic projection has only finitely many double points, called \emph{crossings}, and has no triple points or tangencies. At each crossing, one records which strand passes over and which strand passes under. Some basic examples are shown in Figure~\ref{fig:basic-classical-knots}. 

\begin{figure}[ht]
    \centering
    \includegraphics[width=0.25\textwidth]{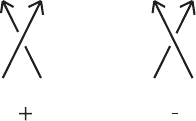}
    \caption{The positive and negative crossings of an oriented knot or link diagram.}
    \label{fig:classical-crossings}
\end{figure}

A knot is called \emph{oriented} if a direction of traversal has been chosen on the embedded circle. More generally, an oriented link is a link in which each component is equipped with such a choice of direction. Once a diagram is oriented, each crossing has a sign, positive or negative, according to the usual convention shown in Figure~\ref{fig:classical-crossings}. Crossing signs are simple local data, but they play an essential role in many classical constructions.

\subsection{Reidemeister moves and diagrammatic equivalence}

The diagrammatic study of knots and links requires choosing a projection to a plane. Such a choice is far from unique: the same embedded knot may be projected onto different planes, and after small perturbations of the embedding or the projection. Consequently, a single knot or link gives rise to many planar diagrams. To use diagrams effectively, one must therefore understand exactly when two diagrams represent ambient isotopic embedded objects. This is the content of the Reidemeister theorem, which translates the geometric problem of ambient isotopy into a combinatorial problem about diagrams and local moves. 

\begin{figure}[ht]
    \centering
    \includegraphics[width=0.72\textwidth]{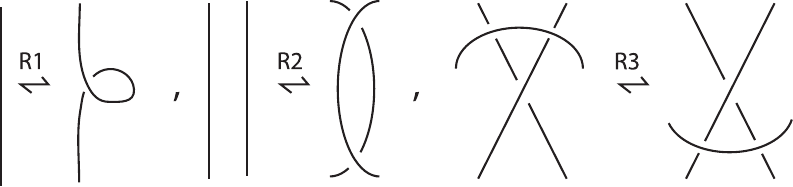}
    \caption{The three Reidemeister moves for classical knot diagrams. These moves generate diagrammatic equivalence for classical knots and links.}
    \label{fig:reidemeister-moves}
\end{figure}

\begin{theorem}[Reidemeister theorem \cite{Reidemeister1932}]\label{ReidThm}
Two classical knot or link diagrams represent ambient isotopic knots or links in \(S^3\) if and only if they are related by a finite sequence of planar isotopies and Reidemeister moves illustrated in Figure~\ref{fig:reidemeister-moves}.
\end{theorem}

The first Reidemeister move, denoted \(R1\), introduces or removes a single twist. The second move, denoted \(R2\), introduces or removes a pair of adjacent crossings, while the third move, denoted \(R3\), slides one strand across a crossing of two other strands. Each move is local, meaning that it changes the diagram only inside a small disk while leaving the rest of the diagram fixed.

In practice, however, Theorem~\ref{ReidThm} is often difficult to apply directly. Given two diagrams, it may be hard to decide whether a sequence of Reidemeister moves relates them, and even harder to find such a sequence explicitly. This difficulty motivates the use of \emph{invariants}. An invariant assigns to each diagram \(D\) an object \(I(D)\) that remains unchanged under all Reidemeister moves. If this condition is satisfied, then \(I(D)\) depends only on the underlying knot or link, not on the particular diagram chosen to represent it.

Invariants provide practical tools for distinguishing knots and links. If two diagrams have different values under some invariant, then they cannot represent ambient isotopic knots or links. The converse is not generally true: equal values of a given invariant do not necessarily imply that the knots or links are equivalent. Thus, invariants are powerful distinguishing tools, although no single standard invariant is sufficient for classification in general.

\begin{example}[Linking number]
Let \(L=L_1\cup L_2\) be an oriented link with two components. The \emph{linking number} of \(L_1\) and \(L_2\) is defined by
\[
\operatorname{lk}(L_1,L_2)
=
\frac{1}{2}\sum_c \operatorname{sign}(c),
\]
where the sum is taken over all crossings between \(L_1\) and \(L_2\), and \(\operatorname{sign}(c)\in\{+1,-1\}\) is the sign of the crossing.

The linking number is invariant under Reidemeister moves and therefore defines an invariant of oriented links.
\end{example}

\begin{example}[Writhe and regular isotopy]
Let \(D\) be an oriented knot or link diagram. The \emph{writhe} of \(D\) is defined by
\[
w(D)=\sum_c \operatorname{sign}(c),
\]
where the sum is taken over all crossings of \(D\). Unlike the linking number, the writhe is not invariant under all Reidemeister moves. A first Reidemeister move changes the writhe by \(+1\) or \(-1\), depending on the sign of the twist introduced or removed. However, the writhe is invariant under the second and third Reidemeister moves, since an \(R2\) move introduces or removes two crossings of opposite sign, while an \(R3\) move preserves the signs of the crossings involved.

Equivalence of diagrams generated by planar isotopy together with the second and third Reidemeister moves is called \emph{regular isotopy} \cite{Kauffman1987,Kauffman1991KnotsPhysics}. Thus, the writhe is an invariant of oriented diagrams under regular isotopy, but not an ambient isotopy invariant of knots or links.
\end{example}

Among the most important and influential examples are polynomial invariants, many of which are defined or characterized by local skein relations. We turn to these next.

\subsection{Polynomial invariants and skein relations}

A polynomial invariant assigns to each knot or link a Laurent polynomial in one or more variables \cite{Alexander1928,Jones1985,FreydYetterHosteLickorishMillettOcneanu1985,
PrzytyckiTraczyk1987,Kauffman1987}. A particularly powerful way to construct polynomial invariants is through \emph{skein relations}. A skein relation is a local relation connecting the values of an invariant on diagrams that are identical outside a small disk and differ only inside that disk. Thus, instead of studying a whole diagram at once, one studies how the invariant changes under controlled local replacements. The Kauffman bracket gives a fundamental example of this idea.

\begin{definition}[Kauffman bracket \cite{Kauffman1987}]
The Kauffman bracket is a map
\[
D\longmapsto \langle D\rangle
\]
from unoriented link diagrams to Laurent polynomials in \(A\), defined recursively by the local relation
\[
\left\langle
\raisebox{-0.15cm}{\includegraphics[height=0.5cm]{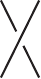}}
\right\rangle
=
A
\left\langle
\raisebox{-0.15cm}{\includegraphics[height=0.5cm]{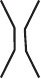}}
\right\rangle
+
A^{-1}
\left\langle
\raisebox{-0.15cm}{\includegraphics[height=0.5cm]{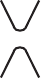}}
\right\rangle,
\]
together with
\[
\langle D\sqcup O\rangle
=
(-A^2-A^{-2})\langle D\rangle
\quad \text{and} \quad
\langle O\rangle=1,
\]
where \(O\) denotes the trivial knot.
\end{definition}

The bracket relation replaces each crossing by one of two smoothings. Repeated application of the relation reduces any diagram to a linear combination of collections of simple closed curves, whose values are determined by the loop relation above \cite{Kauffman1987,Kauffman1991KnotsPhysics}.

The Kauffman bracket is an invariant of regular isotopy, rather than ambient isotopy. This failure is not a defect, but a controlled feature: the change under the first Reidemeister move is exactly compensated by a normalization involving the writhe. In particular, we have that:

\begin{proposition}[\cite{Kauffman1987}]
Let \(D\) be an oriented link diagram representing an oriented link \(L\). Then
\[
f_L(A)=(-A^3)^{-w(D)}\langle D\rangle
\]
is an invariant under all Reidemeister moves and therefore defines an invariant of the oriented link \(L\).
\end{proposition}

\begin{remark}
Under the change of variable \(t=A^{-4}\), the invariant $f_L$ is the Jones polynomial \(V_L(t)\) \cite{Jones1985}, that we will introduce in the next section.
\end{remark}

\begin{figure}[ht]
    \centering
    \includegraphics[width=0.35\textwidth]{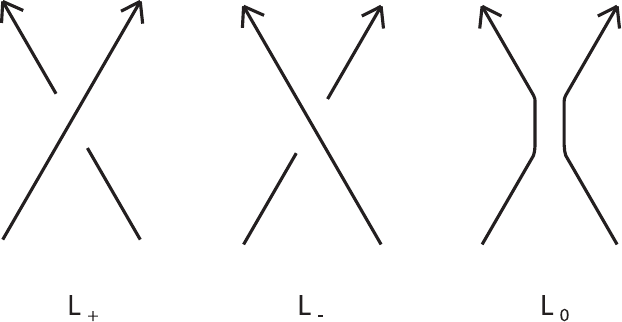}
    \caption{A skein triple consists of three diagrams that differ only inside a small disk.}
    \label{fig:skein-triple}
\end{figure}

The same skein-theoretic philosophy also leads to oriented polynomial invariants. In the oriented setting, one usually considers a skein triple
\[
L_+,\qquad L_-,\qquad L_0,
\]
where the three diagrams are identical outside a small disk and differ inside the disk by a positive crossing, a negative crossing, and an oriented smoothing, respectively (for an illustration see Figure~\ref{fig:skein-triple}).

\begin{theorem}[HOMFLYPT polynomial \cite{FreydYetterHosteLickorishMillettOcneanu1985,PrzytyckiTraczyk1987}]
There exists a unique polynomial invariant \(P_L(\ell,m)\) of oriented links, called the \emph{HOMFLYPT polynomial}, normalized by
\[
P_O(\ell,m)=1,
\]
where \(O\) denotes the unknot, and satisfying the skein relation
\[
\ell P_{L_+}+\ell^{-1}P_{L_-}+mP_{L_0}=0.
\]
\end{theorem}

\begin{remark}
Historically, one of the earliest major polynomial invariants is the Alexander polynomial~\cite{Alexander1928}. For an oriented knot \(K\), it is a Laurent polynomial \(\Delta_K(t)\in \mathbb Z[t,t^{-1}]\), defined up to multiplication by \(\pm t^k\). It admits several equivalent constructions, for example through Seifert matrices or covering-space methods, and remains one of the fundamental classical invariants of knots.

Another important example is the Kauffman polynomial~~\cite{Murakami1987,Kauffman1991KnotsPhysics}, a two-variable invariant of links that also admits a skein-theoretic description. Together with the Alexander, Jones, HOMFLYPT, and Kauffman polynomials \cite{Alexander1928,Jones1985,
FreydYetterHosteLickorishMillettOcneanu1985,
PrzytyckiTraczyk1987,Kauffman1987,Murakami1987}, these invariants form a central family of classical polynomial invariants. They also serve as models for many later constructions in non-classical knot theory, where new crossing types or additional diagrammatic structures require modified skein relations and adapted polynomial invariants.
\end{remark}

The skein-theoretic viewpoint emphasizes local diagrammatic relations. Many of the same polynomial invariants also admit algebraic constructions, in which links are represented as closures of braids and invariance is controlled by Markov moves. We now turn to this braid-theoretic framework.


\subsection{Braids and the Alexander--Markov theorems}

\begin{figure}[ht]
    \centering
    \includegraphics[width=0.5\textwidth]{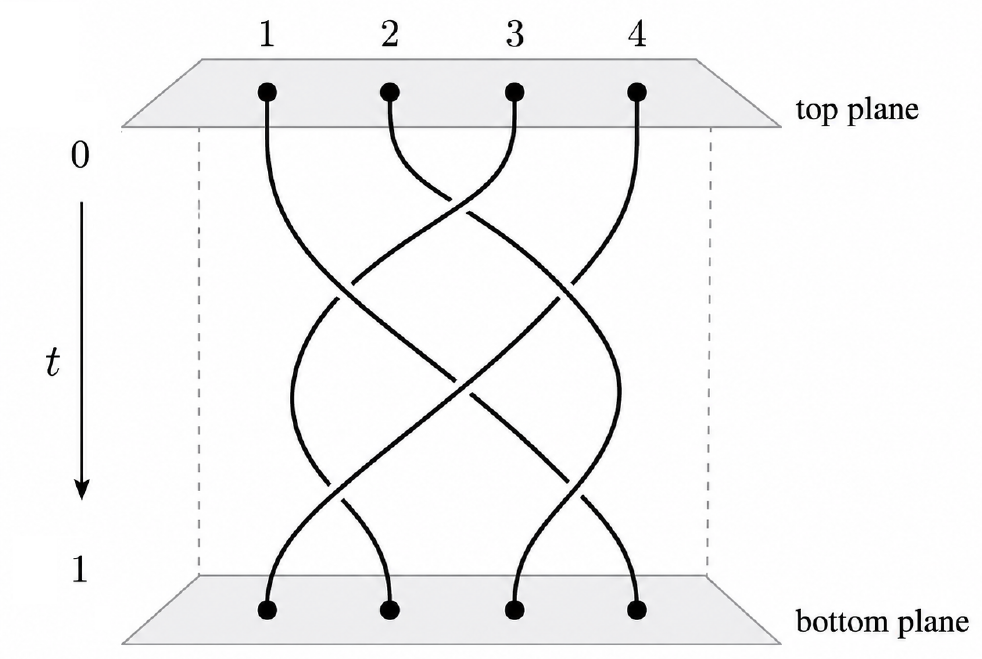}
    \caption{A braid on four strands.}
    \label{fig:braid}
\end{figure}

Braids provide an algebraic way to organize knot and link diagrams. Geometrically, a braid on \(n\) strands may be viewed as a collection of \(n\) disjoint strands in \(\mathbb R^2\times[0,1]\), running monotonically from \(n\) specified points in the top plane to \(n\) specified points in the bottom plane (for an illustration see Figure~\ref{fig:braid}) \cite{Artin1925,Birman1974}. The monotonicity condition is essential: strands may wind around one another and cross in a projection, but they never turn back in the vertical direction.

Two braids are considered equivalent if they are related by an ambient isotopy that preserves the endpoints, and the set of equivalence classes forms the \emph{braid group} \(B_n\). The group operation is given by stacking one braid on top of another, while the inverse of a braid is obtained by reversing its crossings in the appropriate way. This makes braids substantially more algebraic than arbitrary knot diagrams.

\begin{figure}[ht]
    \centering
    \includegraphics[width=0.5\textwidth]{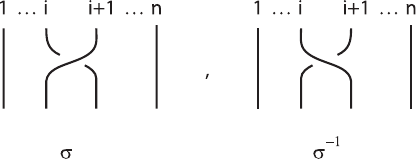}
    \caption{The standard braid generator \(\sigma_i\) and its inverse
    \(\sigma_i^{-1}\).}
    \label{fig:braid-generators}
\end{figure}

The braid group \(B_n\) has standard generators
\[
\sigma_1,\sigma_2,\ldots,\sigma_{n-1},
\]
where \(\sigma_i\) represents the crossing in which the \(i\)-th strand passes over the \((i+1)\)-st strand, and the inverse generator \(\sigma_i^{-1}\) represents the corresponding negative crossing (see Figure~\ref{fig:braid-generators}). 

These generators satisfy two families of relations. First, crossings involving distant pairs of strands commute:
\[
\sigma_i\sigma_j=\sigma_j\sigma_i,
\qquad \text{for } |i-j|>1.
\]
Second, adjacent generators satisfy the braid relation
\[
\sigma_i\sigma_{i+1}\sigma_i
=
\sigma_{i+1}\sigma_i\sigma_{i+1},
\qquad \text{for } 1\leq i\leq n-2.
\]
Together with the existence of inverses \(\sigma_i^{-1}\), these relations give the presentation
\[
B_n=
\left\langle
\sigma_1,\ldots,\sigma_{n-1}
\ \middle|\
\sigma_i\sigma_j=\sigma_j\sigma_i \ (|i-j|>1),\ 
\sigma_i\sigma_{i+1}\sigma_i=
\sigma_{i+1}\sigma_i\sigma_{i+1}
\right\rangle .
\]

\begin{figure}[ht]
    \centering
    \includegraphics[width=0.7\textwidth]{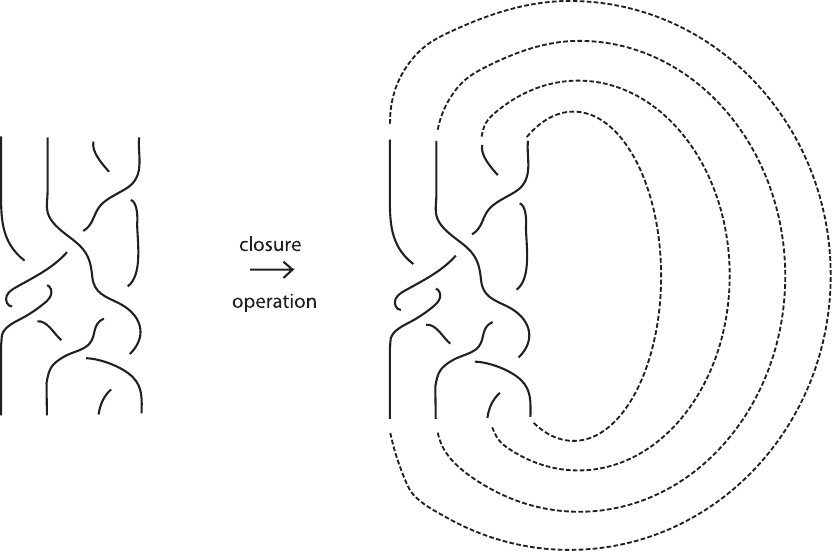}
    \caption{The closure of a braid produces a classical oriented link.}
    \label{fig:braid-closure}
\end{figure}

The connection between braids and links is given by the {\it closure operation}. If \(\beta\in B_n\), its closure, \(\widehat{\beta}\), is obtained by joining corresponding top and bottom endpoints by disjoint arcs, as shown in Figure~\ref{fig:braid-closure}. Since the strands of a braid are naturally oriented from top to bottom, the closure produces an oriented link.

At first sight, braid closures seem to form a restricted class of links, because braid diagrams are required to be monotone in one direction. The Alexander theorem shows that this restriction is not a loss of generality.

\begin{theorem}[Alexander theorem \cite{Alexander1923}]
Every oriented classical link can be represented as the closure of a braid.
\end{theorem}

Alexander's theorem allows one to study oriented links through braid groups. However, the braid representative of a link is far from unique. Different braids, sometimes with different numbers of strands, may have isotopic closures. The precise equivalence relation among braid representatives is
given by Markov's theorem.

\begin{theorem}[Markov theorem \cite{Markov1935}]
Two braids have isotopic closures if and only if they are related by a finite sequence of braid isotopies and the following Markov moves:
\[
\begin{array}{lll}
\text{Conjugation:}
& \beta \sim \alpha\beta\alpha^{-1},
& \alpha,\beta\in B_n, \\[4pt]
\text{Stabilization/destabilization:}
& \beta \sim \beta\sigma_n^{\pm1},
& \beta\in B_n\subset B_{n+1}.
\end{array}
\]
In the second move, \(\beta\) is regarded as an element of \(B_{n+1}\) by adding one trivial strand.
\end{theorem}

\begin{figure}[ht]
    \centering
    \includegraphics[width=0.95\textwidth]{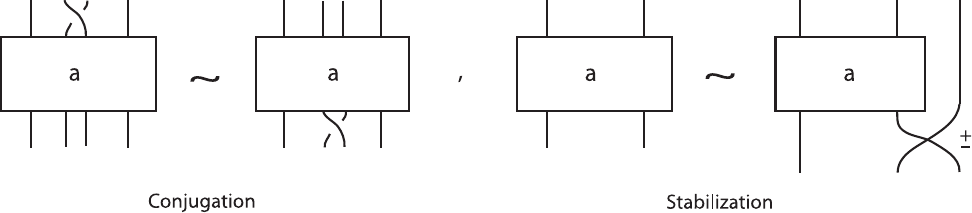}
    \caption{The two types of Markov moves: conjugation and stabilization.}
    \label{fig:markov-moves}
\end{figure}

The two types of Markov moves are illustrated in
Figure~\ref{fig:markov-moves}. Conjugation changes the braid representative without changing the number of strands, while stabilization and destabilization pass between braid groups with different numbers of strands.

Together, the Alexander and Markov theorems form the algebraic counterpart of the Reidemeister theorem. Alexander's theorem says that oriented links may be
represented by braid closures, while Markov's theorem describes exactly when two braid representatives determine the same link. This makes it possible to translate questions about links into questions about braid groups, provided one works modulo Markov equivalence. This braid-theoretic viewpoint is the starting point for many algebraic constructions of link invariants.


\subsection{From braids to link invariants: the algebraic framework}

As mentioned before, the Alexander and Markov theorems translate the construction of oriented link invariants into an algebraic problem on braids \cite{Alexander1923,Markov1935,Birman1974,Jones1987}. A function defined on braids may be an invariant of braids without being an invariant of their closures as links. To descend to a link invariant, it must be unchanged under the Markov moves, namely conjugation and stabilization. Thus, one seeks algebraic constructions on braid groups (or on their quotient algebras) that are compatible with these moves, often after a suitable normalization \cite{Jones1987,FreydYetterHosteLickorishMillettOcneanu1985,
PrzytyckiTraczyk1987}.

A central method is to represent braid groups inside suitable algebras and then apply trace maps \cite{Jones1987,BirmanWenzl1989,Murakami1987}. More precisely, one maps a braid \(\beta\in B_n\) to an element of an algebra \(A_n\), applies a trace-like linear functional \(\operatorname{tr}_n:A_n\to R\) with values in a coefficient ring \(R\), and then normalizes the result so that it is invariant under the Markov moves. If this is done correctly, the resulting element of \(R\) depends only on the isotopy class of the closed braid \(\widehat{\beta}\), and therefore defines an invariant of oriented links.

The basic example is provided by the Hecke algebra of type \(A\). It is obtained from the braid group algebra by imposing a quadratic relation on the braid generators. More precisely, we have:

\begin{definition}[Hecke algebra of type \(A\)]
The Hecke algebra \(H_n(q)\) is the algebra generated by
\[
g_1,\ldots,g_{n-1},
\]
subject to the braid relations
\[
g_i g_j = g_j g_i \qquad \text{for } |i-j|>1,
\]
\[
g_i g_{i+1} g_i = g_{i+1} g_i g_{i+1}
\qquad \text{for } 1\leq i\leq n-2,
\]
and the quadratic relation
\[
g_i^2=(q-1)g_i+q, \quad \text{equivalently} \quad (g_i-q)(g_i+1)=0.
\]
\end{definition}

The defining braid relations allow one to send the braid generator \(\sigma_i\) to \(g_i\). Moreover, the quadratic relation ensures that each \(g_i\) is invertible, with
\[
g_i^{-1}=q^{-1}g_i+(q^{-1}-1).
\]
Thus the braid relations in \(H_n(q)\) allow every braid word in the generators \(\sigma_i^{\pm1}\) to be interpreted as an element of the Hecke algebra by sending
\[
\sigma_i\longmapsto g_i,
\qquad
\sigma_i^{-1}\longmapsto g_i^{-1}.
\]

\begin{remark}
The algebra \(H_n(q)\) is a quotient of the braid group algebra in which each braid generator satisfies a quadratic relation. It may also be viewed as a deformation of the group algebra of the symmetric group: under the specialization \(q=1\), the relation \((g_i-q)(g_i+1)=0\) becomes \((g_i-1)(g_i+1)=0\), which corresponds to the relation \(s_i^2=1\) for the elementary transpositions in the symmetric group.
\end{remark}

To obtain link invariants, one considers the tower of Hecke algebras
\[
H_1(q)\subset H_2(q)\subset H_3(q)\subset\cdots,
\]
where the inclusion sends \(g_i\in H_n(q)\) to the corresponding generator in \(H_{n+1}(q)\). The essential additional ingredient is a Markov trace on this tower.

\begin{definition}[Markov trace \cite{Jones1987}]
A {\it Markov trace} on the tower \(\{H_n(q)\}_{n\geq 1}\) is a family of linear maps
\[
\operatorname{tr}_n:H_n(q)\to R
\]
satisfying, for suitable parameters and for all admissible elements,
\[
\operatorname{tr}_n(ab)=\operatorname{tr}_n(ba),
\]
and if \(a\in A_n\), where \(A_n\) is regarded as a subalgebra of \(A_{n+1}\) under the natural inclusion, then:
\[
\operatorname{tr}_{n+1}(a)=\operatorname{tr}_n(a),
\quad \text{and} \quad
\operatorname{tr}_{n+1}(a g_n)=z\,\operatorname{tr}_n(a).
\]
Usually one also fixes a normalization such as
\(\operatorname{tr}_1(1)=1\).
\end{definition}

The trace condition \(\operatorname{tr}_n(ab)=\operatorname{tr}_n(ba)\) is what makes the construction compatible with conjugation of braids. The Markov property \(\operatorname{tr}_{n+1}(a g_n)=z\,\operatorname{tr}_n(a)\) controls stabilization. Together with the quadratic relation in the Hecke algebra, it also controls negative stabilization. After suitable normalization, the Markov trace on the tower of Hecke algebras yields an invariant of oriented links, namely the HOMFLYPT polynomial.

\begin{theorem}\label{thm:homflypt-trace}
Let \(\alpha\in B_n\) be a braid whose closure is the oriented link \(L=\widehat{\alpha}\). Let
\[
\pi:B_n\longrightarrow H_n(q)
\]
be the canonical homomorphism defined by
\[
\pi(\sigma_i)=g_i.
\]
If \(e=e(\alpha)\) denotes the exponent sum of the braid word \(\alpha\), then the normalized trace expression
\[
X_L(q,\lambda)
=
\left[
-\frac{1-\lambda q}{\sqrt{\lambda}(1-q)}
\right]^{n-1}
\left(\sqrt{\lambda}\right)^e
\operatorname{tr}\bigl(\pi(\alpha)\bigr)
\]
is invariant under the Markov moves and therefore defines an invariant of oriented links in \(S^3\). This invariant is the HOMFLYPT polynomial.
\end{theorem}

The Jones polynomial also fits naturally into this algebraic framework through the Temperley--Lieb algebra \cite{Jones1985,Jones1987,Kauffman1987}. In the Hecke-algebraic approach, the Temperley--Lieb algebra of type \(A\) is obtained as a quotient of the Hecke algebra by an additional local relation.

\begin{definition}[Temperley--Lieb algebra as a Hecke quotient]
The Temperley--Lieb algebra \(TL_n(q)\) is the quotient of \(H_n(q)\) by the two-sided ideal generated by
\[
g_{1,2}
=
1
+q(g_1+g_2)
+q^2(g_1g_2+g_2g_1)
+q^3g_1g_2g_1.
\]
Equivalently,
\[
TL_n(q)=H_n(q)/\langle g_{1,2}\rangle .
\]
\end{definition}

The Ocneanu trace on the Hecke algebras factors through this quotient when the trace annihilates the defining ideal of \(TL_n(q)\). In the normalization used above, this corresponds to imposing
\[
\operatorname{tr}(g_{1,2})=0.
\]
Under the corresponding specialization of the HOMFLYPT-type invariant, one obtains the Jones polynomial, $V_L$. More precisely, with the normalization of Theorem~\ref{thm:homflypt-trace}, setting \(\lambda=q\) gives
\[
V_L(q)=X_L(q,q).
\]

The braid-theoretic construction of polynomial invariants illustrates a broader principle: links may often be studied effectively through braid representatives adapted to the ambient space. For links in other three-manifolds, this leads to generalized braid groups or mixed braid groups, together with algebraic structures reflecting the topology of the manifold. These braid methods provide
an effective approach to skein modules, to which we now turn.



\section{Skein Modules and Knots in Three-Manifolds}
\label{sec:skein-modules}

The preceding section concerned classical links in \(S^3\). A first major generalization is obtained by changing the ambient space itself. Instead of studying links only in \(S^3\), one may consider links embedded in an arbitrary oriented three-manifold \(M\). The local skein relations remain essentially the same, but the topology of the ambient manifold becomes part of the structure \cite{Przytycki1991,Diamantis2025SurveySkeinModules}.

This leads to the theory of skein modules, introduced independently by Przytycki and Turaev~\cite{Przytycki1991,Turaev1990SkeinSolidTorus}. Skein modules extend the idea of polynomial link invariants from \(S^3\) to arbitrary three-manifolds. In \(S^3\), a skein relation often determines a polynomial invariant. In a general three-manifold, the same relation produces an algebraic module
generated by links in that manifold.


\subsection{From polynomial invariants to skein modules}

Let \(M\) be an oriented three-manifold, and let \(\mathcal L(M)\) denote a set of isotopy classes of links in \(M\), possibly with additional structure
such as orientations or framings. A skein module is obtained by taking the free module generated by \(\mathcal L(M)\) and quotienting by local skein
relations. More precisely, we have:

\begin{definition}[Skein module \cite{Przytycki1991}]
Let \(M\) be an oriented three-manifold and let \(R\) be a commutative ring with unit. A \emph{skein module} of \(M\) is an \(R\)-module of the form
\[
R\mathcal L(M)/\mathcal S,
\]
where \(R\mathcal L(M)\) is the free \(R\)-module generated by isotopy classes of links in \(M\), and \(\mathcal S\) is the submodule generated by a chosen
collection of local skein relations.
\end{definition}

The two examples most relevant for this survey are the Kauffman bracket skein module and the HOMFLYPT skein module.

\begin{definition}[Kauffman bracket skein module]
Let \(M\) be an oriented three-manifold and let
\(\mathcal L_{\mathrm{fr}}(M)\) be the set of isotopy classes of unoriented framed links in \(M\). The \emph{Kauffman bracket skein module}
\(\operatorname{KBSM}(M)\) is the module over \(R=\mathbb Z[A^{\pm1}]\) generated by \(\mathcal L_{\mathrm{fr}}(M)\), modulo the relations
\[
L_+ - A L_0 - A^{-1}L_\infty=0
\quad
\text{and}
\quad 
L\sqcup O = (-A^2-A^{-2})L,
\]
where \(O\) denotes the trivially framed unknot in a ball disjoint from \(L\).
\end{definition}

As in the classical case, \(L_+\), \(L_0\), and \(L_\infty\) agree outside a small ball and differ inside the ball as in the usual Kauffman bracket skein relation.

For \(M=S^3\), the Kauffman bracket skein module is one-dimensional and recovers the usual Kauffman bracket polynomial for framed links. In more general manifolds, the same local relation produces a module that may be infinite-dimensional or may contain torsion.

\begin{definition}[HOMFLYPT skein module]
Let \(M\) be an oriented three-manifold and let \(\mathcal L(M)\) be the set of isotopy classes of oriented links in \(M\). The \emph{HOMFLYPT skein module} of \(M\) is the module over \(R=\mathbb Z[u^{\pm1},z^{\pm1}]\) generated by \(\mathcal L(M)\), modulo the oriented skein relation
\[
u^{-1}L_+ - uL_- - zL_0=0.
\]
One usually also allows the empty link and imposes the normalization relation
\[
u^{-1}\emptyset-u\emptyset=zT_1,
\]
where \(T_1\) denotes the trivial knot.
\end{definition}

The Kauffman bracket and HOMFLYPT skein modules are parallel constructions, but they behave differently. The Kauffman bracket theory is unoriented and framed,
while the HOMFLYPT theory is oriented. This distinction becomes important in the braid-theoretic approach, where the HOMFLYPT skein module is naturally
related to Hecke algebras and Markov traces, while the Kauffman bracket skein module is related to Temperley--Lieb quotients.


\subsection{The solid torus and mixed links}

The first nontrivial example is the solid torus \(ST=S^1\times D^2\). It already exhibits features that do not occur in \(S^3\), while still being simple enough to be described explicitly.

\begin{figure}[ht]
\centering
\includegraphics[width=0.85\textwidth]{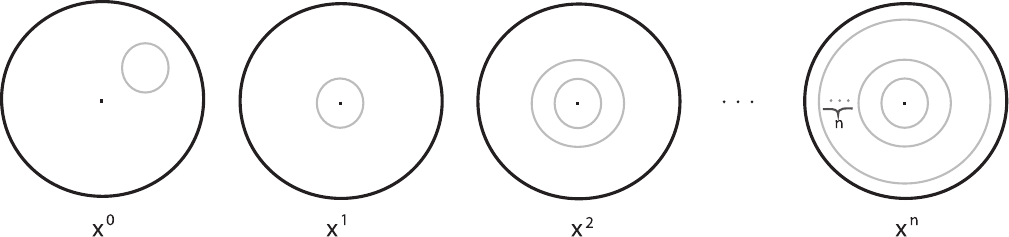}
\caption{A schematic basis for the Kauffman bracket skein module of the solid
torus, generated by parallel copies of the core curve.}
\label{fig:solid-torus-kb-basis}
\end{figure}

The Kauffman bracket skein module of the solid torus is freely generated by parallel copies of the core curve. Thus, a natural basis is given by \(\{x^n\}_{n\geq 0}\), where \(x^n\) denotes \(n\) parallel copies of the core of \(ST\), and \(x^0\) denotes the empty link or affine unknot, depending on the chosen convention (for an illustration see Figure~\ref{fig:solid-torus-kb-basis}).

\begin{figure}[ht]
\centering
\includegraphics[width=0.25\textwidth]{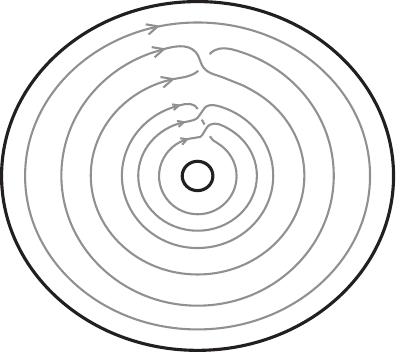}
\caption{A schematic basis element for the HOMFLYPT skein module of the solid torus.}
\label{fig:solid-torus-hom-basis}
\end{figure}

The HOMFLYPT skein module of the solid torus was described by Hoste--Kidwell~\cite{HosteKidwell1990}. It is also freely generated by winding elements around the core of the solid torus, but its oriented nature makes the basis richer than in the Kauffman bracket case (see Figure~\ref{fig:solid-torus-hom-basis}).

\begin{figure}[ht]
\centering
\includegraphics[width=0.2\textwidth]{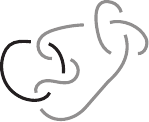}
\caption{A mixed link in \(S^3\) representing a link in the solid torus.}
\label{fig:mixed-link}
\end{figure}

The braid approach to links in \(ST\) uses the fact that the solid torus may be represented inside \(S^3\) by fixing an unknotted component. A link in \(ST\)
is then represented by a mixed link in \(S^3\): one part is fixed and encodes the ambient solid torus, while the other part is the moving link under study (for an illustration see Figure~\ref{fig:mixed-link}).

\begin{definition}[Mixed link \cite{LambropoulouRourke1997Markov3Manifolds}]
A \emph{mixed link} is a link in \(S^3\) decomposed into a fixed sublink and a moving sublink. The fixed part represents the ambient manifold, while the
moving part represents the link inside that manifold.
\end{definition}

\begin{figure}[ht]
\centering
\includegraphics[width=0.55\textwidth]{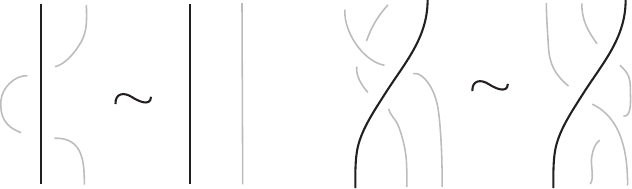}
\caption{Mixed Reidemeister moves for links in the solid torus.}
\label{fig:mixed-reidemeister}
\end{figure}

For the solid torus, the fixed part consists of one unknotted component. The moving part may link with this fixed component, and this linking records how
the link winds around the core of the solid torus. Isotopy in \(ST\) is therefore translated into mixed link isotopy in \(S^3\), that is, ordinary Reidemeister moves on the moving part, together with mixed Reidemeister moves involving both the fixed component and the moving strands, as illustrated in Figure~\ref{fig:mixed-reidemeister} (\cite{LambropoulouRourke1997Markov3Manifolds,
DiamantisLambropoulou2015BraidEquivalence3Manifolds}).

This mixed-link viewpoint is useful because it allows links in \(ST\) to be studied by braid methods \cite{Lambropoulou1999HeckeTypeB}.


\subsection{Mixed braids and type \(B\) algebraic structures}

Passing from mixed links to mixed braids gives the braid-theoretic model for links in the solid torus. A mixed braid consists of one fixed strand, together
with moving strands that braid around it. The corresponding braid group is the Artin braid group of type \(B\), denoted \(B_{1,n}\).

\begin{figure}[ht]
\centering
\includegraphics[width=0.2\textwidth]{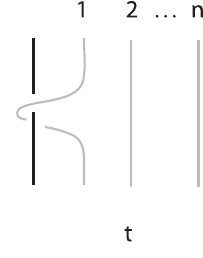}
\caption{The generator \(t\) of the mixed braid group \(B_{1,n}\).}
\label{fig:gent}
\end{figure}

\begin{definition}[Mixed braid group of type \(B\)]
The mixed braid group \(B_{1,n}\) is generated by
\[
t,\sigma_1,\ldots,\sigma_{n-1},
\]
where \(t\) records winding around the fixed strand and the generators \(\sigma_i\) describe the ordinary braiding of the moving strands (for an illustration see Figure~\ref{fig:gent}). It has presentation
\[
B_{1,n} =
\left\langle
t,\sigma_1,\ldots,\sigma_{n-1}
\ \middle|\
\begin{array}{ll}
\sigma_1t\sigma_1t=t\sigma_1t\sigma_1, & \\
t\sigma_i=\sigma_i t, & i>1,\\
\sigma_i\sigma_{i+1}\sigma_i=\sigma_{i+1}\sigma_i\sigma_{i+1},
& 1\leq i\leq n-2,\\
\sigma_i\sigma_j=\sigma_j\sigma_i, & |i-j|>1
\end{array}
\right\rangle .
\]
\end{definition}

\begin{figure}[ht]
\centering
\includegraphics[width=0.55\textwidth]{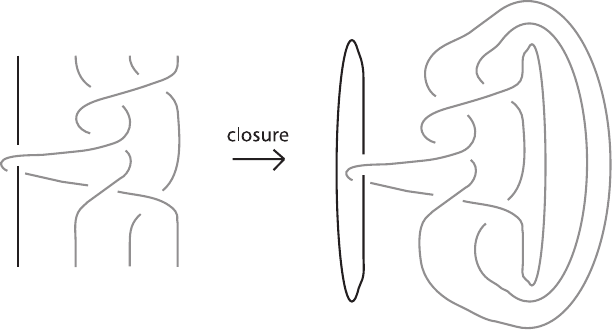}
\caption{The closure of a mixed braid gives a mixed link representing a link in the solid torus.}
\label{fig:mixed-braid-closure}
\end{figure}

The closure of a mixed braid produces a mixed link, and hence a link in the solid torus, as illustrated in Figure~\ref{fig:mixed-braid-closure} (\cite{LambropoulouRourke1997Markov3Manifolds}).

\begin{figure}[ht]
\begin{center}
\includegraphics[width=5.2in]{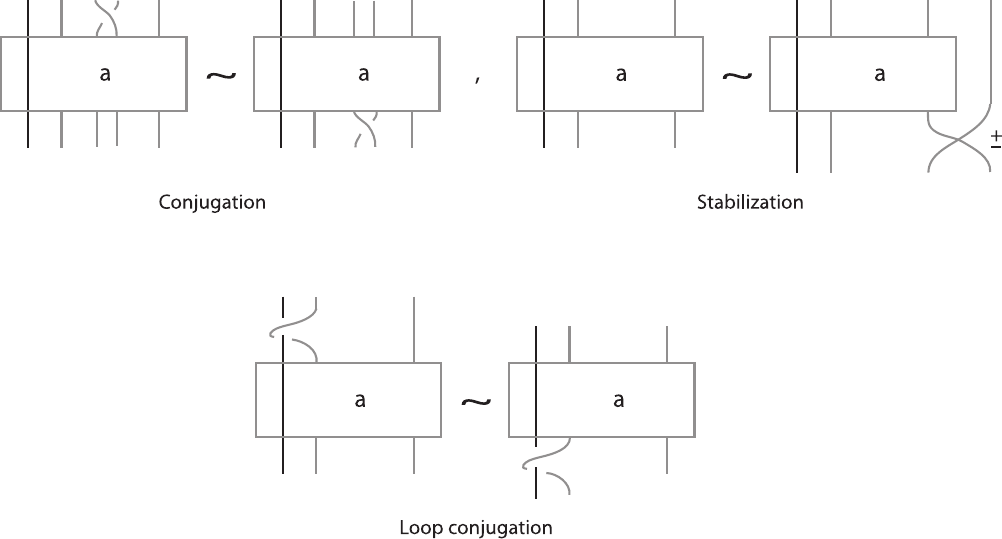}
\end{center}
\caption{Mixed braid equivalence.}
\label{mbeq}
\end{figure}

\begin{theorem}[Alexander and Markov theorems in the solid torus]
Every link in the solid torus may be represented as the closure of a mixed braid. Moreover, two mixed braids represent isotopic links in \(ST\) if and
only if they are related by braid isotopy, conjugation, stabilization, and loop conjugation moves~\cite{LambropoulouRourke1997Markov3Manifolds,
DiamantisLambropoulou2015BraidEquivalence3Manifolds} (for an illustration see Figure~\ref{mbeq}).
\end{theorem}

Thus, the classical braid description of links in \(S^3\) extends to links in the solid torus by replacing the ordinary braid group \(B_n\) with the type \(B\) mixed braid group \(B_{1,n}\). The corresponding algebraic object is the generalized Hecke algebra of type \(B\) (\cite{Lambropoulou1999HeckeTypeB}).

\begin{definition}[Generalized Hecke algebra of type \(B\)]
The generalized Hecke algebra \(H_{1,n}(q)\) is obtained from the group algebra of \(B_{1,n}\) by imposing the quadratic relations
\[
g_i^2=(q-1)g_i+q
\]
on the generators corresponding to the moving braid generators \(\sigma_i\). The looping generator \(t\) satisfies no polynomial relation.
\end{definition}

The absence of a polynomial relation for \(t\) makes \(H_{1,n}(q)\) infinite-dimensional. This reflects the topology of the solid torus: links may wind around the core arbitrarily many times.

\begin{figure}[ht]
\centering
\includegraphics[width=0.5\textwidth]{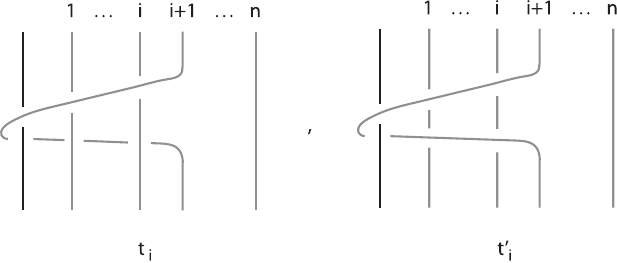}
\caption{Looping elements in the type \(B\) braid setting.}
\label{fig:looping-elements}
\end{figure}

To describe the corresponding trace construction, one introduces looping elements
\[
t_0=t,
\qquad
t_i' = g_i\cdots g_1\,t\,g_1^{-1}\cdots g_i^{-1},
\]
which encode loops around the fixed strand, and which are illustrated in Figure~\ref{fig:looping-elements}.

A Markov trace on the tower of generalized Hecke algebras \(H_{1,n}(q)\) was constructed by Lambropoulou~\cite{Lambropoulou1999HeckeTypeB}. Besides cyclicity and the usual stabilization rule, it includes trace parameters \(s_k\) that record the looping information:
\[
\operatorname{tr}(ab)=\operatorname{tr}(ba),
\qquad
\operatorname{tr}(ag_n)=z\operatorname{tr}(a),
\qquad
\operatorname{tr}\left(a{t_n'}^k\right)=s_k\operatorname{tr}(a).
\]

\begin{theorem}[HOMFLYPT-type invariant of \(ST\) \cite{Lambropoulou1999HeckeTypeB}]
Let \(\mathcal L\) denote the set of oriented links in the solid torus, and let
\(\alpha\in B_{1,n}\) be a mixed braid whose closure is the oriented link
\(\widehat{\alpha}\subset ST\). Let
\[
\pi:B_{1,n}\longrightarrow H_{1,n}(q)
\]
be the canonical map defined by
\[
\pi(t)=t,
\qquad
\pi(\sigma_i)=g_i.
\]
If \(e=e(\alpha)\) denotes the exponent sum of the braid generators
\(\sigma_i\) in \(\alpha\), then
\[
X_{\widehat{\alpha}}
=
\Delta^{\,n-1}
\left(\sqrt{\lambda}\right)^e
\operatorname{tr}\bigl(\pi(\alpha)\bigr),
\]
where
\[
\Delta
=
-\frac{1-\lambda q}{\sqrt{\lambda}(1-q)},
\qquad
\lambda
=
\frac{z+1-q}{qz}.
\]
This defines an invariant of oriented links in \(ST\). 
\end{theorem}

This theorem makes explicit how the Markov trace on the type \(B\) Hecke algebras produces a universal invariant for oriented links in the solid torus:
the trace evaluates the mixed braid representative, while the normalization corrects the behavior under Markov stabilization.

\begin{remark}\rm 
The universal HOMFLYPT-type invariant $X$ recovers the HOMFLYPT skein module of the solid torus.    
\end{remark}

The Kauffman bracket side has an analogous algebraic formulation \cite{Diamantis2024KBSMLensUnoriented, Diamantis2019AlternativeBasisST}. Following the classical passage from Hecke algebras to Temperley--Lieb algebras, one defines a generalized Temperley--Lieb algebra of type \(B\), denoted \(TL_{1,n}\), as a quotient of \(H_{1,n}(q)\) by the appropriate Temperley--Lieb ideal. Under suitable restrictions on the trace parameters, the Markov trace on \(H_{1,n}(q)\) factors through this quotient.

More precisely, in the presentation with parameters \(u\) and \(v\), the trace factors through \(TL_{1,n}\) when
\[
z=-\frac{1}{u(1+u^2)}.
\]
For this value of \(z\), one has \(\lambda=u^4\), and the following invariant is obtained.

\begin{theorem}[Kauffman bracket type invariant of \(ST\) \cite{Diamantis2024KBSMLensUnoriented}]
Let \(\alpha\in B_{1,n}\) be a mixed braid whose closure is the link \(\widehat{\alpha}\subset ST\). Let
\[
\overline{\pi}:B_{1,n}\longrightarrow TL_{1,n}
\]
be the canonical map defined by
\[
\overline{\pi}(t)=t,
\qquad
\overline{\pi}(\sigma_i)=g_i.
\]
If \(e=e(\alpha)\) denotes the exponent sum of the braid generators \(\sigma_i\) in \(\alpha\), then
\[
V_{\widehat{\alpha}}(u,v)
=
\left(-\frac{1+u^2}{u}\right)^{n-1}
u^{2e}
\operatorname{tr}\bigl(\overline{\pi}(\alpha)\bigr)
\]
defines a universal Kauffman bracket type invariant for links in the solid torus.
\end{theorem}

This invariant recovers the Kauffman bracket skein module of \(ST\): it separates the elements of the standard basis of \(\operatorname{KBSM}(ST)\). In particular, for the basis elements represented by
\(t t_1'\cdots t_n'\), one has
\[
\operatorname{tr}(t t_1'\cdots t_n')=s_1^n.
\]

\begin{remark}
The solid torus mirrors the classical situation. In \(S^3\), type \(A\) braid groups lead to Hecke and Temperley--Lieb algebras, and hence to the HOMFLYPT
and Kauffman bracket invariants. In the solid torus, type \(B\) mixed braid groups lead to type \(B\) Hecke and Temperley--Lieb structures, and hence to the HOMFLYPT and Kauffman bracket skein modules of \(ST\) (\cite{Jones1987,Lambropoulou1999HeckeTypeB,
DiamantisLambropoulou2016NewBasisHOMFLYPTST,
Diamantis2019AlternativeBasisST}).
\end{remark}


\subsection{Lens spaces and braid band moves}

The solid torus is also the starting point for studying skein modules of other three-manifolds. By the Lickorish--Wallace theorem, every closed, connected, orientable three-manifold may be obtained by surgery on a link in \(S^3\) \cite{Lickorish1997}. Consequently, links in any three-manifold may be studied through mixed links and mixed braids, with additional moves encoding the surgery description \cite{LambropoulouRourke1997Markov3Manifolds,
DiamantisLambropoulou2015BraidEquivalence3Manifolds}.

\begin{figure}[ht]
\begin{center}\includegraphics[width=3.1in]{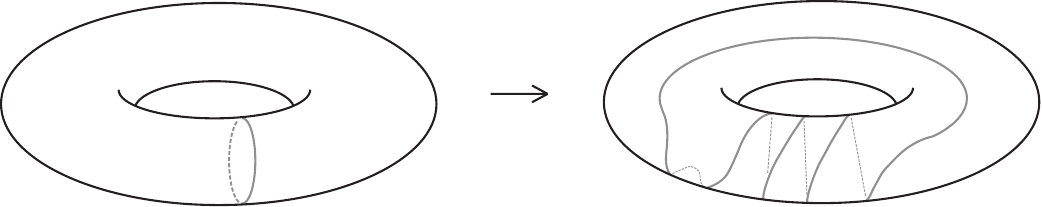}
\end{center}
\caption{The homeomorphism $h$.}
\label{st1}
\end{figure}

Lens spaces provide the simplest family of examples. The lens space \(L(p,1)\) is obtained from \(S^3\) by integral surgery on the unknot with surgery coefficient \(p\). Equivalently, one removes a solid torus and glues it back so that a meridian is identified with a \((p,1)\)-curve on the boundary of the complementary solid torus (see Figure~\ref{st1}).

\begin{figure}[ht]
\begin{center}
\includegraphics[width=2.4in]{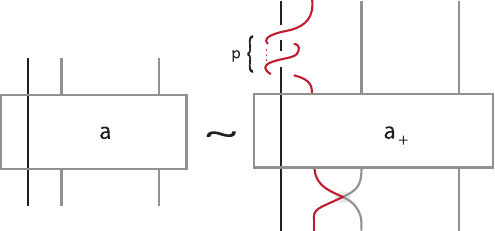}
\end{center}
\caption{Braid band moves.}
\label{bbm}
\end{figure}

In mixed link language, isotopy in \(L(p,1)\) is described by isotopy in the solid torus together with band moves reflecting the surgery. On the braid
level, these become \emph{braid band moves}, and they are illustrated in Figure~\ref{bbm}. Thus, the Markov theorem for links in \(L(p,1)\) is obtained from the Markov theorem for mixed braids in the solid torus by adding braid band moves to the list of allowed equivalence moves~\cite{LambropoulouRourke1997Markov3Manifolds,
DiamantisLambropoulou2015BraidEquivalence3Manifolds}.

Consequently, skein modules of lens spaces can be approached by starting with skein modules of the solid torus and imposing additional relations induced by
braid band moves (\cite{DiamantisLambropoulouPrzytycki2016TopologicalSteps,
DiamantisLambropoulou2017BraidApproach,
DiamantisLambropoulou2019ImportantStep,
Diamantis2021HOMFLYPTSubmodulesLens,
Diamantis2024KBSMLensUnoriented}). Schematically, one passes from the skein module of \(ST\) to that of \(L(p,1)\) by adding relations of the form
\[
a=bbm(a),
\]
where \(a\) ranges over suitable basis elements and \(bbm(a)\) denotes the result of performing a braid band move.

For the Kauffman bracket skein module this method gives an explicit finite basis.

\begin{theorem}[\cite{Diamantis2024KBSMLensUnoriented}]
For \(p\neq 0\), the Kauffman bracket skein module of \(L(p,1)\) admits the natural basis
\[
\mathcal B_p=\{t^n\mid 0\leq n\leq \lfloor p/2\rfloor\}.
\]
\end{theorem}

The HOMFLYPT skein module of \(L(p,1)\) can be approached in the same spirit, but the computation is more subtle (\cite{Diamantis2015PhD,
DiamantisLambropoulouPrzytycki2016TopologicalSteps,
DiamantisLambropoulou2017BraidApproach,
DiamantisLambropoulou2019ImportantStep,
Diamantis2021HOMFLYPTSubmodulesLens}). One starts from the HOMFLYPT skein
module of the solid torus and imposes the additional relations induced by braid band moves. Equivalently, one seeks to normalize the universal HOMFLYPT-type invariant \(X\) of \(ST\) so that it is invariant under all braid band moves.

For this purpose, it is useful to work with the basis
\[
\Lambda=
\left\{
t^{k_0}t_1^{k_1}\cdots t_n^{k_n}
\ \middle|\
n\in\mathbb N,\;
k_i\in\mathbb Z\setminus\{0\},\;
k_i\geq k_{i+1}
\right\},
\]
whose elements have no gaps in the indices, ordered exponents, and no braiding tails~\cite{DiamantisLambropoulou2016NewBasisHOMFLYPTST}. This basis is adapted to braid band moves and replaces earlier bases in which gaps, unordered exponents, and braiding tails make the relations harder to control.

The computation of \(\operatorname{HOM}(L(p,1))\) then reduces to solving the infinite system of equations
\[
X_{\widehat{\tau}}
=
X_{\widehat{bbm_i(\tau)}},
\qquad
\tau\in\Lambda,
\]
where \(bbm_i(\tau)\) denotes the result of performing a braid band move on the \(i\)-th moving strand of \(\tau\). In schematic form, one may write
\[
\operatorname{HOM}(L(p,1))
\sim
\frac{\Lambda}
{\langle \tau-bbm_i(\tau)\rangle},
\qquad
\tau\in\Lambda,
\]
with the understanding that this notation represents the infinite system of relations imposed on the universal invariant \(X\), rather than a simple finite
presentation. Thus, the braid approach reduces the HOMFLYPT skein module of \(L(p,1)\) to an algebraic problem: determine the consequences of all braid band move relations on the basis \(\Lambda\) (\cite{Diamantis2015PhD,DiamantisLambropoulou2019ImportantStep,
Diamantis2021HOMFLYPTSubmodulesLens}).

This illustrates a recurring theme in the braid approach to skein modules: topological information about the ambient manifold is encoded by additional
geometric moves, while skein-theoretic information is encoded algebraically through bases, Hecke-type algebras, and trace parameters (\cite{Diamantis2025SurveySkeinModules}).


\subsection{Torsion phenomena and \(S^1\times S^2\)}

The manifold \(S^1\times S^2\), which is the lens space \(L(0,1)\), is an important test case for skein-module methods because its skein modules exhibit torsion phenomena (\cite{Diamantis2024KBSMS1S2,Diamantis2025HOMFLYPTS1S2}). This shows that skein modules contain information about the ambient manifold that has no direct analogue in the classical polynomial setting of \(S^3\).

In the braid-theoretic approach, links in \(S^1\times S^2\) are represented by mixed braids in the solid torus. The passage from \(ST\) to \(S^1\times S^2\)
is encoded by braid band moves corresponding to the \(0\)-surgery description. Thus, extending an invariant from \(ST\) to \(S^1\times S^2\) requires imposing
invariance under these braid band moves.

This strategy has been carried out for both the Kauffman bracket and the HOMFLYPT skein modules of \(S^1\times S^2\) (\cite{Diamantis2024KBSMS1S2,Diamantis2025HOMFLYPTS1S2}). In the Kauffman bracket case, one
starts from the universal Kauffman bracket type invariant of the solid torus and imposes braid band move relations. The resulting computation shows that
\[
\operatorname{KBSM}(S^1\times S^2)/\operatorname{Tor}
\cong
\mathbb Z[A^{\pm1}],
\]
so the free part is generated by the empty link, or equivalently by the unknot. The torsion is detected through the core element \(t\) of the solid torus; in particular, the computation shows that the nontrivial powers of the core become torsion after imposing the \(0\)-surgery relations~\cite{Diamantis2024KBSMS1S2}.

The HOMFLYPT skein module of \(S^1\times S^2\) can be treated in an analogous way. Starting from the universal HOMFLYPT-type invariant \(X\) of the solid
torus, one imposes the infinite system of braid band move equations \( X_{\widehat{\alpha}} = X_{\widehat{bbm(\alpha)}}\), where \(\alpha\) ranges over a basis of the HOMFLYPT skein module of \(ST\). Solving this system yields a decomposition of the form
\[
\operatorname{HOM}(S^1\times S^2)
\cong
R\{\emptyset\}\oplus \mathcal T,
\]
where \(R\{\emptyset\}\) is the rank-one free submodule generated by the empty link and \(\mathcal T\) is a torsion submodule~\cite{Diamantis2025HOMFLYPTS1S2}. Thus, as in the Kauffman bracket case, the free part is generated by the empty link, while the remaining information is torsion.

The examples of the solid torus and lens spaces show how the classical braid--algebra--trace framework extends beyond \(S^3\) \cite{Diamantis2025SurveySkeinModules}. The ambient space is encoded by fixed components, mixed braid groups, and additional moves such as braid band moves. Thus, changing the ambient manifold provides a first major route from classical knot theory to generalized knot theories.


\section{Pseudo and Singular Knot Theories}

A second major way to extend classical knot theory is to modify the local structure of crossings. In a classical diagram, every crossing carries complete
over/under information. Pseudo and singular knot theories both depart from this classical setting by introducing new types of local crossing data, but they do so with different interpretations (\cite{Hanaki2010,HenrichEtAl2013Pseudoknots,
Vassiliev1990,Baez1992,Gemein1997SingularBraids}).

Pseudo knot theory replaces some classical crossings by unresolved crossings. Such crossings record ambiguity or incomplete information: locally, the diagram does not specify which strand passes over and which strand passes under. Singular knot theory, on the other hand, replaces some classical crossings by
prescribed transverse double points. These are not ambiguous crossings, but local singularities that are treated as part of the structure of the diagram.

Thus, pseudo and singular knot theories are conceptually distinct. Pseudo knot theory studies controlled ambiguity, while singular knot theory studies controlled degeneration. Nevertheless, the two theories are algebraically close. Both admit similar equivalence moves, braid-theoretic formulations, Alexander and Markov type theorems, and Hecke-type algebraic structures (\cite{Gemein1997SingularBraids,FennKeymanRourke1998SingularBraid,
ParisRabenda2004,BardakovJablanWang2016PseudoBraids,
Diamantis2026HOMFLYPTPseudoResolution}). This makes them a natural pair to study together.


\subsection{Pseudo and singular crossings}

We begin by defining the two local objects that distinguish pseudo and singular diagrams from classical knot diagrams.

\begin{definition}[Pseudo crossing and singular crossing] 
A \emph{pseudo crossing} (or pre-crossing) is a double point of a diagram at which no over/under information is specified. A \emph{singular crossing} is a prescribed transverse double point, treated as part of the geometric structure of the diagram. \end{definition}

A pseudo crossing represents a place where one has not chosen between the positive and negative classical crossings, while a singular crossing represents a genuine double point in an immersed curve, and this double point remains part of the object under consideration. This distinction is important throughout the theory, since the two theories lead to different interpretations of local replacement rules, skein relations, and polynomial invariants. 

\begin{figure}[ht]
\centering
\includegraphics[width=0.15\textwidth]{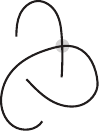}
\caption{A pseudo knot.}
\label{fig:pseudo-knot}
\end{figure}

Pseudo and singular knot theories arose from different motivations. Pseudo knot theory was introduced by Hanaki~\cite{Hanaki2010} and developed further by Dye, Henrich, Hoberg, Jablan, Kauffman, and others
\cite{Dye2010,HenrichEtAl2013Pseudoknots,BardakovJablanWang2016PseudoBraids}, as a way to study diagrams in which some crossing information is intentionally left unresolved. Singular knot theory has an earlier origin in the study of immersed curves and finite-type invariants, especially through the work of Vassiliev, Birman, Baez, Bar-Natan, and others  \cite{Vassiliev1990,Baez1992,BarNatan1995,Gemein1997SingularBraids,
FennKeymanRourke1998SingularBraid,ParisRabenda2004}. Although the two theories have different geometric interpretations, they are naturally compared because both replace the classical crossing by a new local object and both lead to generalized Reidemeister moves, braid monoids, skein-type relations, and polynomial invariants.


\subsection{Diagrams and generalized Reidemeister moves}

A \emph{pseudo knot diagram} is a classical knot diagram in which some crossings are replaced by pseudo crossings. A \emph{pseudo link diagram} is defined
analogously. Similarly, a \emph{singular knot diagram} is a classical knot diagram containing ordinary classical crossings together with finitely many singular crossings; singular link diagrams are defined in the same way.

Every classical knot diagram may be regarded as a pseudo diagram with no pseudo crossings, and also as a singular diagram with no singular crossings. At the
opposite extreme, a diagram all of whose crossings are pseudo crossings is closely related to a \emph{knot shadow}, namely a generic immersion of a circle
in the plane without over/under information at its double points.

\begin{figure}[ht]
    \centering
    \includegraphics[width=0.99\textwidth]{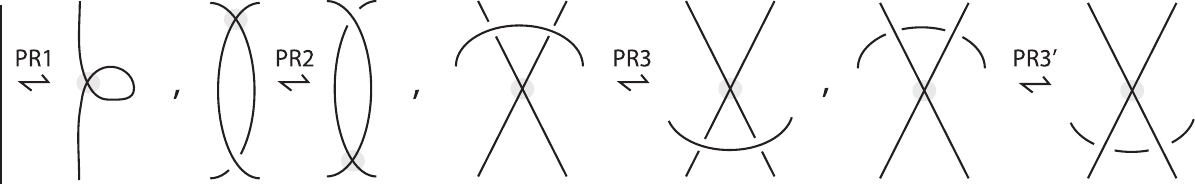}
    \caption{The pseudo Reidemeister moves. These moves extend the classical
    Reidemeister moves by allowing local transformations involving pseudo
    crossings.}
    \label{fig:pseudo-reidemeister}
\end{figure}

The equivalence relations in the two theories are obtained by extending the classical Reidemeister moves, but the extensions reflect the different meanings of the new crossings. For pseudo knots, one allows the classical Reidemeister moves together with pseudo Reidemeister moves involving pseudo crossings as illustrated in Figure~\ref{fig:pseudo-reidemeister}. These moves express the fact that pseudo crossings behave like unresolved classical crossings rather than fixed singular vertices.

For singular knots, one allows the classical Reidemeister moves together with singular Reidemeister moves involving singular crossings. These moves describe how prescribed double points interact with ordinary strands and ordinary crossings. Unlike the pseudo case, there is no singular analogue of the first
pseudo Reidemeister move PR1: a singular crossing cannot simply be introduced or removed as an unresolved kink. This reflects the fact that singular crossings are part of the geometric data of the diagram.

\begin{definition}[Pseudo and singular knots]
A \emph{pseudo knot} is an equivalence class of pseudo knot diagrams under classical and pseudo Reidemeister moves. A \emph{singular knot} is an equivalence class of singular knot diagrams under classical and singular
Reidemeister moves.
\end{definition}

Thus, pseudo diagrams interpolate between classical knot diagrams and knot shadows, whereas singular diagrams interpolate between embedded curves and
immersed curves with prescribed transverse double points, and this distinction is already visible at the level of moves.


\subsection{Resolution and tangle insertion viewpoints}

The difference between pseudo and singular crossings becomes especially clear when one tries to relate the generalized diagrams back to classical knot theory. In both cases, one can replace a generalized crossing by classical local data. However, the meaning of this replacement is different.

For pseudo knots, the natural operation is resolution \cite{HenrichEtAl2013Pseudoknots}. Each pseudo crossing may be replaced by either a positive or a negative classical crossing. Thus, a pseudo diagram represents a family of classical diagrams obtained by resolving all its pseudo crossings.

\begin{definition}[Resolution states and resolution set \cite{HenrichEtAl2013Pseudoknots}]
Let \(D\) be a pseudo knot or link diagram with \(m\) pseudo crossings. A \emph{classical resolution state} of \(D\) is obtained by replacing each pseudo
crossing by one of the two classical crossings. The collection of all such states is denoted by \(\mathcal S(D)\).

The set of classical knot or link types represented by the diagrams in \(\mathcal S(D)\) is called the \emph{resolution set} of \(D\) and is denoted
by \(\mathcal R(D)\).
\end{definition}

A diagram with \(m\) pseudo crossings has \(2^m\) classical resolution states, although different states may represent the same classical knot or link type.
This leads naturally to a weighted version of the resolution set.

\begin{figure}[ht]
    \centering
    \includegraphics[width=0.45\textwidth]{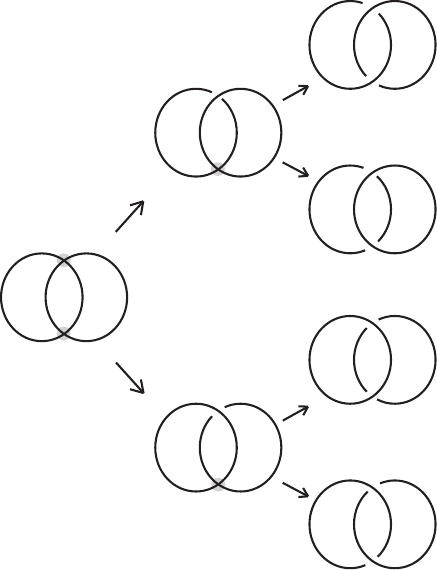}
    \caption{Resolution tree for a pseudo link diagram with two pseudo crossings. The four classical links obtained by resolving the two pseudo
    crossings consist of Hopf links and unlinks.}
    \label{fig:weighted-resolution}
\end{figure}

\begin{definition}[Weighted resolution set \cite{HenrichEtAl2013Pseudoknots}]
Let \(K_1,\ldots,K_r\) be the distinct classical knot or link types appearing among the \(2^m\) resolution states of \(D\). If \(n_i\) resolution states represent \(K_i\), then the \emph{weighted resolution set}, or \emph{WeRe set}, of \(D\) is
\[
\operatorname{WeRe}(D)
=
\{(K_1,n_1/2^m),\ldots,(K_r,n_r/2^m)\}.
\]
\end{definition}

The WeRe set, introduced by Henrich et al. \cite{HenrichEtAl2013Pseudoknots}, refines the ordinary
resolution set by recording multiplicities. It remembers not only which classical knot types appear, but also how often they occur among all classical
resolutions. For example, in Figure~\ref{fig:weighted-resolution}, two of the four resolution states yield the Hopf link, while the remaining two yield the unlink. Hence
\[
\operatorname{WeRe}(D)
=
\{(\text{Hopf link},1/2),(\text{unlink},1/2)\}.
\]

Resolution sets also provide a direct way to extend classical invariants to pseudo knots. If \(I\) is an invariant of classical knots or links, then one
may associate to a pseudo diagram \(D\) the multiset
\[
\{I(D')\mid D'\in \mathcal S(D)\}.
\]
Equivalently, if
\(\operatorname{WeRe}(D)=\{(K_1,p_1),\ldots,(K_r,p_r)\}\), and \(P_K\) is a polynomial invariant of a classical knot or link \(K\), one may define the averaged invariant
\[
P_D^{\mathrm{avg}}
=
\sum_{i=1}^r p_i P_{K_i}.
\]
This construction reflects the interpretation of a pseudo crossing as unresolved classical crossing information.

A broader formulation is given by the tangle insertion method of Henrich and Kauffman \cite{HenrichKauffman2017TangleInsertion}. In this approach, a generalized crossing is replaced by a chosen tangle, or by a formal combination of tangles, and classical invariants are then applied after insertion. For pseudo knots, such insertions encode possible local resolutions or weighted choices.

For singular knots, the corresponding replacement is not probabilistic or ambiguous, but signed (\cite{Vassiliev1990,BarNatan1995}). In the theory of finite-type invariants, a classical
invariant \(I\) is extended to singular diagrams by the singular skein relation
\[
I(L_\times)=I(L_+)-I(L_-),
\]
where \(L_\times\), \(L_+\), and \(L_-\) denote three link diagrams that are identical outside a small disk and differ inside the disk by a singular crossing, a positive crossing, and a negative crossing, respectively. Thus, a singular crossing is interpreted as the difference between the two classical crossings. This signed expansion is the algebraic expression of the fact that singular crossings model controlled degeneration rather than ambiguity.

\begin{definition}[Finite-type invariant]
Let \(I\) be an invariant of classical knots, extended to singular knots by \(I(D_\times)=I(D_+)-I(D_-)\). We say that \(I\) is a \emph{finite-type invariant of type at most \(k\)} if this extension vanishes on every singular knot with more than \(k\) singular
crossings.
\end{definition}

Thus, singular crossings provide a way to organize classical knot invariants according to how they behave under controlled singular degenerations.


\subsection{Braid and algebraic structures}

Pseudo and singular knot theories also admit braid-theoretic formulations \cite{Gemein1997SingularBraids,FennKeymanRourke1998SingularBraid,
BardakovJablanWang2016PseudoBraids, Diamantis2021TiedPseudoLinks}. As in the classical setting, the goal is to represent diagrammatic objects by
braid-like objects and then describe equivalence of closures through generalized Markov moves. The resulting braid objects are monoids rather than
groups, because the generators corresponding to pseudo or singular crossings are not invertible.

\begin{definition}[Pseudo braid monoid \cite{BardakovJablanWang2016PseudoBraids}]
The \emph{pseudo braid monoid} \(PM_n\) is generated by
\[
\sigma_1^{\pm1},\ldots,\sigma_{n-1}^{\pm1},
\qquad
p_1,\ldots,p_{n-1},
\]
where \(\sigma_i^{\pm1}\) are the classical braid generators and \(p_i\) represents a pseudo crossing between the \(i\)-th and \((i+1)\)-st strands. The generators \(\sigma_i\) satisfy the usual braid relations and the pseudo generators satisfy the mixed relations
\[
\begin{array}{rcll}
p_i\, p_j & = & p_j\, p_i, & {\rm if}\ |i-j|\geq 2\\
&&&\\
p_i\, \sigma_j^{\pm 1} & = & \sigma_j^{\pm 1}\, p_i, & {\rm if}\ |i-j|\geq 2\\
&&&\\
p_i\, \sigma_i^{\pm 1} & = & \sigma_i^{\pm 1}\, p_i, & i=1, \ldots, n-1\\
&&&\\
\sigma_i\, \sigma_{i+1}\, p_i & = & p_{i+1}\, \sigma_i\, \sigma_{i+1}, & i=1, \ldots, n-2\\
&&&\\
\sigma_{i+1}\, \sigma_i\, p_{i+1} & = & p_{i}\, \sigma_{i+1}\, \sigma_i, & i=1, \ldots, n-2\\
\end{array}
\]
\end{definition}

\begin{definition}[Singular braid monoid \cite{Birman1974}]
The \emph{singular braid monoid} \(SM_n\) is generated by
\[
\sigma_1^{\pm1},\ldots,\sigma_{n-1}^{\pm1},
\qquad
\tau_1,\ldots,\tau_{n-1},
\]
where the generators \(\sigma_i^{\pm1}\) are the classical braid generators and \(\tau_i\) represents a singular crossing between the \(i\)-th and \((i+1)\)-st strands. The defining relations are obtained from the presentation of \(PM_n\) above by replacing each pseudo generator \(p_i\) with the singular generator \(\tau_i\).
\end{definition}

\begin{figure}[ht]
    \centering
    \includegraphics[width=2.4in]{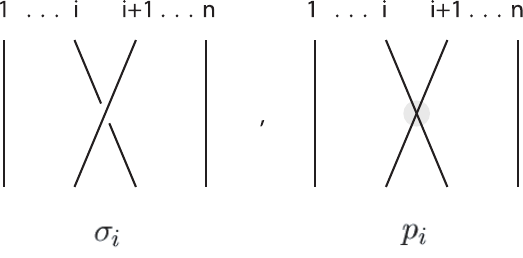}
    \caption{Generalized braid generators. Besides the classical braid generator \(\sigma_i\), pseudo braid theory introduces \(p_i\), representing
    an unresolved crossing, while singular braid theory introduces \(\tau_i\), representing a prescribed transverse double point.}
    \label{fig:generalized-braid-generators}
\end{figure}

Although pseudo and singular crossings have different diagrammatic meanings, their braid monoids are algebraically very close.

\begin{theorem}
The pseudo braid monoid \(PM_n\) is isomorphic to the singular braid monoid \(SM_n\). The isomorphism is given on generators by
\[
\sigma_i \longmapsto \sigma_i,
\qquad
p_i \longmapsto \tau_i .
\]
\end{theorem}

Note that at the monoid level, pseudo and singular braid theories therefore have the same formal algebraic structure. The interpretation, however, remains different.

In \cite{FennKeymanRourke1998SingularBraid} it is shown that $SM_n$ embeds in a group, the singular braid group $SB_n$. It follows that $PM_n$ embeds in a group also, the pseudo braid group $PB_n$. 

\begin{theorem}
The singular braid monoid embeds in a group. Consequently, via the isomorphism \(PM_n\cong SM_n\), the pseudo braid monoid also embeds in a group.
\end{theorem}

The closure of a generalized braid is defined as in the classical case by joining corresponding top and bottom endpoints. The closure of a pseudo braid
is a pseudo link, while the closure of a singular braid is a singular link \cite{Diamantis2021TiedPseudoLinks}.

\begin{theorem}[Alexander-type theorems]
Every oriented pseudo link is the closure of a pseudo braid. Similarly, every singular link is the closure of a singular braid.
\end{theorem}

We now recall the Markov-type theorems for singular and pseudo links in \(S^3\). These results extend the classical Markov theorem by adding moves that account for the non-classical generators \(\tau_i\) and \(p_i\)
\cite{BardakovJablanWang2016PseudoBraids, Diamantis2021TiedPseudoLinks}.

Let
\[
SM_\infty=\bigcup_{n\geq 1}SM_n,
\qquad
PM_\infty=\bigcup_{n\geq 1}PM_n,
\]
where the inclusions are given by adding a trivial strand.

\begin{theorem}[Markov theorem for singular links]\label{thm:singular-markov}
Let \(\alpha\in SM_n\) and \(\beta\in SM_m\) be singular braids. Then the closures \(\widehat{\alpha}\) and \(\widehat{\beta}\) are isotopic as oriented singular links if and only if \(\alpha\) and \(\beta\) are related in \[ SM_\infty=\bigcup_{n\geq 1}SM_n \] by singular braid relations and a finite sequence of the following moves:

\[ \begin{array}{rcll}
\text{Commuting:} & \alpha_1\alpha_2 &\sim& \alpha_2\alpha_1, \qquad \alpha_1,\alpha_2\in SM_n, \\[0.3cm]
\text{Stabilization:} & \alpha &\sim& \alpha\sigma_n^{\pm1}, \qquad \alpha\in SM_n. \end{array}
\]

Here \(\alpha\sigma_n^{\pm1}\) is regarded as an element of \(SM_{n+1}\) after adding one trivial strand.
\end{theorem}

\begin{theorem}[Markov theorem for pseudo links]\label{thm:pseudo-markov}
Two pseudo braids have isotopic closures as pseudo links if and only if one can be obtained from the other by a finite sequence of pseudo braid relations
and the following moves:
\[
\begin{array}{rcll}
\text{Conjugation:}
&
\alpha
&\sim&
\beta^{\pm1}\alpha\beta^{\mp1},
\qquad
\alpha\in PM_n,\ \beta\in B_n,
\\[0.3cm]
\text{Commuting:}
&
\alpha\beta
&\sim&
\beta\alpha,
\qquad
\alpha,\beta\in PM_n,
\\[0.3cm]
\text{Stabilization:}
&
\alpha
&\sim&
\alpha\sigma_n^{\pm1},
\qquad
\alpha\in PM_n,
\\[0.3cm]
\text{Pseudo-stabilization:}
&
\alpha
&\sim&
\alpha p_n,
\qquad
\alpha\in PM_n.
\end{array}
\]
Here \(\alpha\sigma_n^{\pm1}\) and \(\alpha p_n\) are regarded as elements of \(PM_{n+1}\) after adding one trivial strand.
\end{theorem}

\begin{remark}
The Markov theorem for pseudo knots also admits an equivalent \(L\)-move formulation (\cite{Lambropoulou2007LMoves, Diamantis2021TiedPseudoLinks}): the conjugation and stabilization moves may be replaced by \(L\)-moves, while one retains the commuting move and the pseudo-stabilization move. This formulation is often more flexible when extending braid equivalence to other diagrammatic settings.
\end{remark}

These results provide the entry point for algebraic constructions, since generalized braid monoids may be mapped to Hecke-type quotient algebras and combined with trace maps to produce polynomial invariants \cite{ParisRabenda2004,Diamantis2026HOMFLYPTPseudoResolution}.


\subsection{Hecke-type algebras and trace constructions}

The analogues of the Alexander and Markov theorems allow pseudo and singular links to be studied through braid representatives. As in the classical setting, this opens the door to algebraic constructions based on quotient algebras and suitable trace functionals. The main idea is to map generalized braid monoids to
Hecke-type algebras and then normalize the resulting trace expressions so that they are invariant under the corresponding Markov moves.

Singular braid monoids lead to singular Hecke algebras, introduced and studied by Paris and Rabenda~\cite{ParisRabenda2004}. These algebras extend the classical Hecke algebra by adjoining generators corresponding to singular crossings and imposing relations compatible with the singular braid monoid.
They are naturally connected with signed skein expansions and finite-type invariants.

Pseudo braid monoids lead instead to pseudo Hecke-type algebras \cite{Diamantis2023PseudoSingularSolidTorus}. In this setting, the additional generators encode unresolved crossing information and may be interpreted through resolution maps or tangle insertion \cite{Diamantis2026HOMFLYPTPseudoResolution}.

\begin{definition}[Pseudo Hecke algebra of type \(A\) \cite{Diamantis2023PseudoSingularSolidTorus}]
The \emph{pseudo Hecke algebra of type \(A\)}, denoted
\(\mathcal{PH}_n(q)\), is generated by
\[
g_1,\ldots,g_{n-1},
\qquad
p_1,\ldots,p_{n-1},
\]
where the generators \(g_i\) satisfy the classical Hecke relations
\[
g_i g_j=g_j g_i \qquad (|i-j|\geq 2),
\]
\[
g_i g_{i+1}g_i=g_{i+1}g_i g_{i+1}
\qquad (i=1,\ldots,n-2),
\]
and
\[
g_i^2=(q-1)g_i+q.
\]
The generators \(p_i\) correspond to pseudo crossings and satisfy the relations induced by the pseudo braid monoid, with \(\sigma_i\) replaced by \(g_i\).
\end{definition}

Thus, the canonical map \(PM_n\longrightarrow \mathcal{PH}_n(q)\) sends \(\sigma_i\longmapsto g_i\) and \(p_i\longmapsto p_i\). The Hecke relation makes each \(g_i\) invertible, with \(g_i^{-1}=q^{-1}g_i+(q^{-1}-1)\).

One algebraic way to incorporate the unresolved nature of pseudo crossings is through a resolution homomorphism. In recent work by the author \cite{Diamantis2026HOMFLYPTPseudoResolution}, one defines
\[
\rho_{X,Y}:\mathcal{PH}_n(q)\longrightarrow H_n(q)
\]
by
\[
\rho_{X,Y}(g_i)=g_i,
\qquad
\rho_{X,Y}(p_i)=Xg_i+Yg_i^{-1},
\]
where \(X\) and \(Y\) are parameters.

This map is the algebraic counterpart of the resolution viewpoint. A pseudo crossing is not treated as a singular vertex; rather, it is replaced by a
linear combination of the two corresponding classical Hecke generators. If
\[
\operatorname{tr}_n:H_n(q)\to R
\]
denotes the Ocneanu trace, we set
\[
T_n=\operatorname{tr}_n\circ \rho_{X,Y}.
\]
The family \(\{T_n\}_{n\geq 1}\) behaves like a Markov trace adapted to pseudo crossings.

\begin{proposition}[\cite{Diamantis2026HOMFLYPTPseudoResolution}]
For all \(a,b\in \mathcal{PH}_n(q)\), the maps \(T_n\) satisfy
\[
T_n(ab)=T_n(ba).
\]
Moreover, under the natural inclusion \(\mathcal{PH}_n(q)\subset \mathcal{PH}_{n+1}(q)\), one has
\[
T_{n+1}(a)=T_n(a),
\]
\[
T_{n+1}(a g_n)=z\,T_n(a),
\]
and
\[
T_{n+1}(a p_n)=(Xz+Yz_-)\,T_n(a),
\]
where
\[
z_-:=q^{-1}z+q^{-1}-1.
\]
\end{proposition}

The last relation records the trace-theoretic contribution of a terminal pseudo crossing. Indeed, under the resolution homomorphism,
\[
p_n\longmapsto Xg_n+Yg_n^{-1},
\]
so a pseudo stabilization decomposes into a linear combination of positive and negative classical stabilizations. The coefficient \(Xz+Yz_-\) is precisely the factor that must be compensated in the normalization of the pseudo link invariant.

Let \(\alpha\in PM_n\). We denote by \(e(\alpha)\) the exponent sum of the classical braid generators \(\sigma_i\) in a word representative of \(\alpha\),
and by \(d(\alpha)\) the number of pseudo generators \(p_i\). Both quantities are well-defined on \(PM_n\), since the defining relations preserve the classical exponent sum and the number of pseudo generators.

We now define
\[
\mathcal P(\widehat{\alpha})
=
A^{n-1}B^{e(\alpha)}C^{d(\alpha)}T_n(\alpha),
\]
where \(T_n(\alpha)\) means \(T_n\) applied to the image of \(\alpha\) in \(\mathcal{PH}_n(q)\). The constants \(A\), \(B\), and \(C\) are chosen so that
\(\mathcal P\) is invariant under the pseudo Markov moves.

The stabilization conditions are
\[
ABz=1,
\qquad
AB^{-1}z_-=1,
\qquad
AC(Xz+Yz_-)=1.
\]
Equivalently,
\[
B^2=\frac{z_-}{z},
\qquad
A=\frac{1}{Bz},
\qquad
C=\frac{1}{A(Xz+Yz_-)}.
\]
We work over a coefficient ring containing the relevant parameters and inverses, and containing a choice of \(B\) satisfying
\[
B^2=\frac{z_-}{z}.
\]

\begin{theorem}[HOMFLYPT-type invariant of pseudo links \cite{Diamantis2026HOMFLYPTPseudoResolution}]
Let \(\mathcal P\) be the function defined above by
\[
\mathcal P(\widehat{\alpha})
=
A^{n-1}B^{e(\alpha)}C^{d(\alpha)}T_n(\alpha),
\qquad \alpha\in PM_n,
\]
where \(T_n=\operatorname{tr}_n\circ\rho_{X,Y}\). If the constants \(A,B,C\) satisfy the normalization conditions above, then \(\mathcal P\) is invariant
under the pseudo Markov moves and therefore defines an invariant of oriented pseudo links.
\end{theorem}

The proof follows the same pattern as in the classical case. The trace property gives invariance under conjugation and commuting moves, while the three
normalization conditions give invariance under positive stabilization, negative stabilization, and pseudo-stabilization. Thus the expression is invariant under all pseudo Markov moves.

This invariant extends the classical HOMFLYPT construction. The new feature is the pseudo degree \(d(\alpha)\), whose normalization factor \(C^{d(\alpha)}\) reflects the role of the first pseudo Reidemeister move, or equivalently the pseudo-stabilization move in the Markov theorem.

The construction may be summarized schematically as
\[
PM_n
\longrightarrow
\mathcal{PH}_n(q)
\xrightarrow{\ \rho_{X,Y}\ }
H_n(q)
\xrightarrow{\ \operatorname{tr}\ }
R,
\]
followed by the normalization \(A^{n-1}B^{e(\alpha)}C^{d(\alpha)}\).

The invariant also satisfies a natural pseudo skein relation. Let \(L_p\), \(L_+\), and \(L_-\) be three oriented pseudo link diagrams that are identical outside a small disk and differ inside the disk by a pseudo crossing, a positive crossing, and a negative crossing, respectively. Then
\[
\mathcal P(L_p)
=
\lambda_+\,\mathcal P(L_+)
+
\lambda_-\,\mathcal P(L_-),
\]
where
\[
\lambda_+=XCB^{-1},
\qquad
\lambda_-=YCB.
\]
This relation is the diagrammatic expression of the algebraic resolution
\[
p_i\longmapsto Xg_i+Yg_i^{-1}.
\]

\begin{table}[ht]
\centering
\[
\begin{array}{c|c|c}
 & \text{Pseudo theory} & \text{Singular theory} \\
\hline
\text{extra crossing} & \text{unresolved crossing} & \text{prescribed double point} \\
\text{braid generator} & p_i & \tau_i \\
\text{braid object} & \text{pseudo braid monoid} & \text{singular braid monoid} \\
\text{replacement principle} & \text{choices/weights} & \text{signed difference} \\
\text{algebraic behavior} & p_i \mapsto Xg_i+Yg_i^{-1} & \text{singular/finite-type expansion}
\end{array}
\]
\caption{Comparison of the pseudo and singular braid--algebra frameworks.}
\label{tab:pseudo-singular-comparison}
\end{table}

Together, these constructions show how pseudo and singular knot theories extend the classical braid--algebra--trace mechanism in two different directions.
Pseudo theory algebraizes ambiguity through resolution maps and pseudo normalization, while singular theory algebraizes controlled degeneration through singular skein expansions and finite-type invariants.


\subsection{Pseudo and singular structures in the solid torus}

\begin{figure}[ht]
\centering
\includegraphics[width=0.18\textwidth]{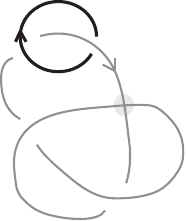}
\caption{A mixed pseudo link in \(S^3\). Replacing the pseudo crossings by singular crossings gives the corresponding mixed singular setting.}
\label{fig:mixed-pseudo-link}
\end{figure}

The braid-theoretic framework for links in the solid torus extends naturally to pseudo and singular knot theory \cite{Diamantis2023PseudoSingularSolidTorus}. We view \(ST\) as the complement of a solid torus in \(S^3\) and a pseudo link in \(ST\) as a
\emph{mixed pseudo link} in \(S^3\) (for an illustration see Figure~\ref{fig:mixed-pseudo-link}). Similarly, singular links in \(ST\) may be represented as mixed singular links in \(S^3\).

Note that we do not allow special crossings between the fixed and moving parts: pseudo crossings, in the pseudo case, and singular crossings, in the singular case, occur only among the moving strands.

The Reidemeister theorem for links in \(ST\) is therefore replaced by a mixed version.

\begin{theorem}[Reidemeister theorem for mixed pseudo links \cite{Diamantis2023PseudoSingularSolidTorus}]
Two mixed pseudo links in \(S^3\) represent isotopic pseudo links in \(ST\) if and only if their diagrams differ by a finite sequence of classical Reidemeister moves and pseudo Reidemeister moves on the moving part, together with mixed Reidemeister moves involving the fixed component and the moving part.
\end{theorem}

\begin{figure}[ht]
\centering
\includegraphics[width=0.28\textwidth]{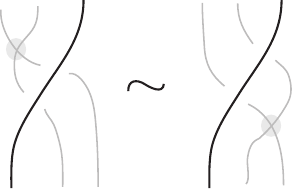}
\caption{Mixed Reidemeister moves involving the fixed component and the moving part.}
\label{fig:mixed-pseudo-reidemeister}
\end{figure}

For singular links in \(ST\), one obtains the analogous statement by replacing pseudo crossings by singular crossings and by using the singular Reidemeister moves. In particular, the pseudo \(PR1\)-move is not present in the singular setting.

\medskip

\subsubsection*{Mixed pseudo and singular braids}

A mixed pseudo braid consists of a fixed strand \(I\), representing the complementary solid torus, together with a moving pseudo braid on \(n\) strands. After applying the usual parting procedure, one may assume that the fixed strand is the first strand and that the remaining \(n\) strands form the moving part. The closure is defined by closing the moving braid strands in the usual way while keeping the fixed strand fixed.

\begin{definition}[Mixed pseudo braid]
A \emph{mixed pseudo braid} on \(n\) moving strands is a braid of the form \(I\cup \beta\), where \(I\) is the fixed identity strand and \(\beta\) is a pseudo braid on the moving strands. Its closure is a mixed pseudo link in \(S^3\), and hence a pseudo link in \(ST\).
\end{definition}

\begin{figure}[ht]
\centering
\includegraphics[width=0.6\textwidth]{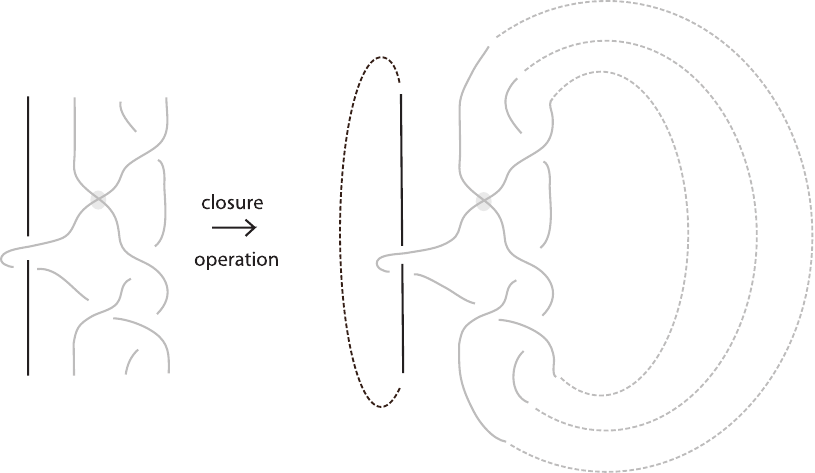}
\caption{The closure of a mixed pseudo braid to a mixed pseudo link.}
\label{fig:mixed-pseudo-braid-closure}
\end{figure}

The braiding algorithm for mixed pseudo and mixed singular links follows the classical braiding algorithm for mixed links, with one additional precaution: special crossings must be arranged so that the braiding process does not alter them. More precisely, pseudo or singular crossings involving up-arcs are first rotated so that the incident arcs are directed downward; then the usual braiding algorithm is applied to the remaining up-arcs, keeping the fixed
component unchanged.

This gives Alexander-type theorems for pseudo and singular links in \(ST\).

\begin{theorem}[Alexander theorem for pseudo and singular links in \(ST\) \cite{Diamantis2023PseudoSingularSolidTorus}]
Every oriented pseudo link in \(ST\) is isotopic to the closure of a mixed pseudo braid. Similarly, every oriented singular link in \(ST\) is isotopic to the closure of a mixed singular braid.
\end{theorem}

\medskip

\subsubsection*{The mixed pseudo braid monoid}

The algebraic counterpart of mixed pseudo braids is the mixed pseudo braid monoid of type \(B\). It extends the type \(B\) braid group \(B_{1,n}\) by
adding pseudo generators among the moving strands.

\begin{definition}[Mixed pseudo braid monoid]
The \emph{mixed pseudo braid monoid of type \(B\)}, denoted \(PM_{1,n}\), is generated by
\[
t^{\pm1},\quad
\sigma_1^{\pm1},\ldots,\sigma_{n-1}^{\pm1},
\quad
p_1,\ldots,p_{n-1}.
\]
Here \(t\) is the looping generator around the fixed strand, the \(\sigma_i^{\pm1}\) are the classical braid generators on the moving strands, and \(p_i\) represents a pseudo crossing between the \(i\)-th and
\((i+1)\)-st moving strands.

The generators \(\sigma_i^{\pm1}\) and \(p_i\) satisfy the relations of the
pseudo braid monoid \(PM_n\). In addition, the loop generator satisfies the
type \(B\) relations
\[
t\sigma_i=\sigma_i t,
\quad \text{and} \quad
tp_i=p_i t,
\quad \text{for}\ i>1,
\]
\[
t\sigma_1t\sigma_1=\sigma_1t\sigma_1t,
\quad \text{and} \quad 
t\sigma_1tp_1=p_1t\sigma_1t.
\]
\end{definition}

The mixed singular braid monoid is obtained from the same formal presentation by replacing pseudo generators by singular generators.

\begin{definition}[Mixed singular braid monoid]
The \emph{mixed singular braid monoid of type \(B\)}, denoted \(SM_{1,n}\), is generated by
\[
t^{\pm1},\quad
\sigma_1^{\pm1},\ldots,\sigma_{n-1}^{\pm1},
\quad
\tau_1,\ldots,\tau_{n-1},
\]
with defining relations obtained from those of \(PM_{1,n}\) by replacing each pseudo generator \(p_i\) by a singular generator \(\tau_i\).
\end{definition}

As in the case of \(S^3\), the pseudo and singular mixed braid monoids are formally isomorphic, although their generators have different geometric interpretations.

\begin{theorem}[\cite{Diamantis2023PseudoSingularSolidTorus}]
There is an isomorphism
\[
SM_{1,n}\longrightarrow PM_{1,n}
\]
defined by
\[
\sigma_i^{\pm1}\longmapsto \sigma_i^{\pm1},
\qquad
t^{\pm1}\longmapsto t^{\pm1},
\qquad
\tau_i\longmapsto p_i.
\]
\end{theorem}

Thus, at the level of monoid presentations, the singular and pseudo type \(B\) braid theories have the same formal algebraic structure.

\begin{remark}
The solid torus is the first example of a broader family of ambient three-manifold settings in which pseudo links can be studied by mixed-link and mixed-braid methods. In particular, pseudo links in handlebodies can be treated by replacing the single fixed component of the solid-torus picture with an
appropriate fixed sublink encoding the handlebody. This leads to corresponding mixed pseudo braid descriptions and generalized Reidemeister-type moves for pseudo links in handlebodies~\cite{Diamantis2021PseudoHandlebodies}.
\end{remark}

\medskip

\subsubsection*{Markov-type theorems in the solid torus}

The Markov theorem for mixed pseudo braids is obtained by adapting the Markov theorem for pseudo braids to the relative setting in which the fixed strand is
kept unchanged (\cite{BardakovJablanWang2016PseudoBraids,
Diamantis2023PseudoSingularSolidTorus}).

\begin{theorem}[Markov theorem for mixed pseudo braids \cite{Diamantis2023PseudoSingularSolidTorus}]
Two mixed pseudo braids have equivalent closures as pseudo links in \(ST\) if and only if one can be obtained from the other by a finite sequence of the
following moves:
\[
\begin{array}{rcll}
\text{Commuting:}
&
\alpha p_i
&\sim&
p_i\alpha,
\qquad
\alpha\in PM_{1,n},
\\[0.3cm]
\text{Conjugation:}
&
\beta
&\sim&
\gamma^{\pm1}\beta\gamma^{\mp1},
\qquad
\beta\in PM_{1,n},\ \gamma\in B_{1,n},
\\[0.3cm]
\text{Real stabilization:}
&
\alpha
&\sim&
\alpha\sigma_n^{\pm1},
\qquad
\alpha\in PM_{1,n},
\\[0.3cm]
\text{Pseudo-stabilization:}
&
\alpha
&\sim&
\alpha p_n,
\qquad
\alpha\in PM_{1,n}.
\end{array}
\]
The stabilization moves are regarded in \(PM_{1,n+1}\) after adding one moving strand.
\end{theorem}

For mixed singular braids, the corresponding theorem is obtained by replacing \(p_i\) with \(\tau_i\) and omitting pseudo-stabilization.

\begin{theorem}[Markov theorem for mixed singular braids \cite{Diamantis2023PseudoSingularSolidTorus}]
Two mixed singular braids have equivalent closures as singular links in \(ST\) if and only if one can be obtained from the other by a finite sequence of the
following moves:
\[
\begin{array}{rcll}
\text{Commuting:}
&
\alpha \tau_i
&\sim&
\tau_i\alpha,
\qquad
\alpha\in SM_{1,n},
\\[0.3cm]
\text{Conjugation:}
&
\beta
&\sim&
\gamma^{\pm1}\beta\gamma^{\mp1},
\qquad
\beta\in SM_{1,n},\ \gamma\in B_{1,n},
\\[0.3cm]
\text{Real stabilization:}
&
\alpha
&\sim&
\alpha\sigma_n^{\pm1},
\qquad
\alpha\in SM_{1,n}.
\end{array}
\]
Again, stabilization is regarded after adding one moving strand.
\end{theorem}

\medskip 

\subsubsection*{Hecke-type algebras and pseudo bracket invariants in \(ST\)}

On the algebraic side, one obtains pseudo and singular analogues of the type \(B\) Hecke algebras (\cite{Lambropoulou1999HeckeTypeB,
Diamantis2023PseudoSingularSolidTorus}). The generalized pseudo Hecke algebra of type \(B\) is
defined as
\[
\mathcal{PH}_{1,n}(q)
=
R[PM_{1,n}]
\Big/
\left\langle
g_i^2-(q-1)g_i-q
\right\rangle,
\]
where \(g_i\) denotes the image of \(\sigma_i\) and the loop generator \(t\) satisfies no polynomial relation in the generalized type \(B\) setting.

The pseudo Hecke algebra of type \(B\) maps naturally to the generalized Hecke algebra \(H_{1,n}(q)\) by the resolution homomorphism
\[
\rho_{X,Y}:\mathcal{PH}_{1,n}(q)\longrightarrow H_{1,n}(q),
\]
defined by
\[
\rho_{X,Y}(g_i)=g_i,
\qquad
\rho_{X,Y}(p_i)=Xg_i+Yg_i^{-1},
\qquad
\rho_{X,Y}(t)=t.
\]
The loop generator \(t\) preserves the winding information around the core of the solid torus, while each pseudo generator is resolved as a linear combination of the positive and negative Hecke generators (\cite{Diamantis2026HOMFLYPTPseudoResolution}).

This provides the natural type \(B\) analogue of the pseudo Hecke construction in \(S^3\). To obtain link invariants, one must combine this resolution map with a Markov trace on \(H_{1,n}(q)\) and a normalization compatible with the mixed pseudo Markov moves. In this way, pseudo links in \(ST\) connect the crossing-based generalizations of this section with the mixed-braid and skein-module framework of Section~\ref{sec:skein-modules}.

\begin{remark}
As in the classical pseudo Kauffman bracket construction of Dye~\cite{Dye2010}, there is also a diagrammatic Kauffman bracket type invariant for pseudo links in \(ST\)~\cite{Diamantis2023PseudoSingularSolidTorus}. Starting from the pseudo bracket polynomial for pseudo links in \(S^3\), one adds the solid-torus relation recording the loop around the core. After writhe normalization, one obtains an invariant
\[
P_K(A,V,s)
=
(-A^{-3})^{w(K)}\langle K\rangle
\]
of oriented pseudo knots in \(ST\), where \(w(K)\) is the writhe of the classical crossings and \(\langle K\rangle\) is the pseudo bracket polynomial. This invariant extends the normalized pseudo bracket polynomial from \(S^3\) to the solid torus.
\end{remark}


\section{Stuck Knots}

The pseudo and singular theories discussed in the previous section modify classical crossings by replacing the usual over/under data with unresolved or
singular local data. Stuck knot theory takes a different route. A stuck crossing is still a fully specified classical crossing, with definite over/under information, but it carries an additional constraint that restricts the way it may participate in isotopy moves (\cite{Bataineh2020,CenicerosElhamdadiKomissarLahrani2024,
CenicerosElhamdadiMagillRosario2023,Diamantis2026StuckKnots}).

Thus, stuck crossings do not encode ambiguity, as pseudo crossings do, and they do not represent prescribed transverse double points, as singular
crossings do. Instead, they retain the classical crossing data while adding a local rigidity condition. In this sense, stuck knot theory is a crossing-based
generalization in which local constraints become part of the topological structure.


\subsection{Stuck crossings and restricted isotopy}

\begin{definition}[Stuck crossing and stuck diagram \cite{Bataineh2020}]
A \emph{stuck crossing} is a classical crossing equipped with an additional marking indicating that the crossing is constrained. A \emph{stuck knot
diagram} is a knot diagram containing ordinary classical crossings together with finitely many stuck crossings.
\end{definition}

\begin{figure}[ht]
\centering
\includegraphics[width=0.7\textwidth]{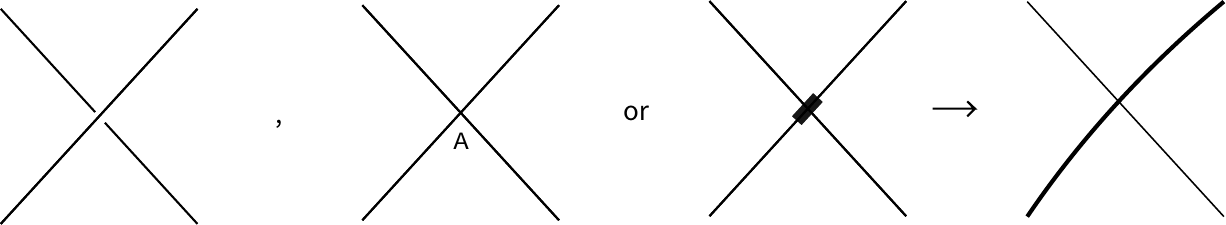}
\caption{A classical crossing and a stuck crossing represented diagrammatically as a rigid-height vertex together with a distinguished region (marked \(A\)) determining the height convention.}
\label{fig:crossing_types}
\end{figure}

\begin{figure}[ht]
\centering
\includegraphics[width=0.95\textwidth]{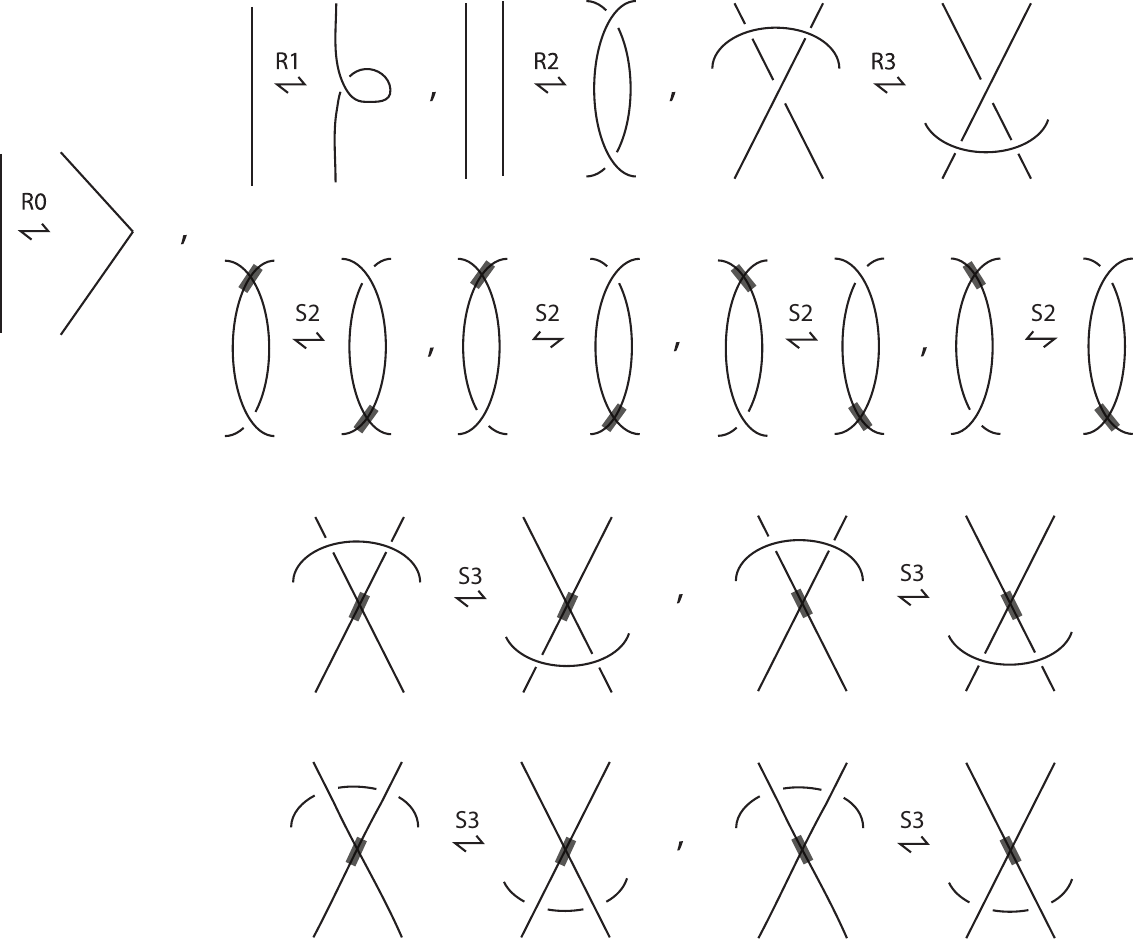}
\caption{Isotopy moves for classical and stuck crossings.}
\label{fig:rigid_moves}
\end{figure}

The equivalence relation for stuck diagrams is generated by the usual Reidemeister moves away from stuck crossings, together with local moves involving stuck crossings that preserve the constrained structure; see Figure~\ref{fig:rigid_moves}. Moves that would remove a stuck crossing or freely change its constrained role are not allowed, unless one explicitly passes to a relaxed theory in which the stuck crossing is released \cite{Diamantis2026StuckKnots}.

\begin{definition}[Stuck knot]
A \emph{stuck knot} is an equivalence class of stuck knot diagrams under the restricted isotopy generated by the moves illustrated in Figure~\ref{fig:rigid_moves}.
\end{definition}

A distinctive feature of stuck knot theory is that two diagrams may have the same underlying classical knot type after the stuck structure is forgotten, while still representing different stuck knots. The difference is not the classical isotopy type of the underlying diagram, but the placement and behavior of the constrained crossings.

\begin{remark}
Pseudo crossings encode unresolved information, singular crossings encode controlled degeneration, and stuck crossings encode constrained classical crossing data. In all three theories, the local crossing structure is no longer governed solely by the classical Reidemeister moves.
\end{remark}


\subsection{A HOMFLYPT-type invariant}

The constructions and results in this subsection and the next one follow \cite{Diamantis2026StuckKnots}. We recall them here in order to illustrate how rigidity constraints can be incorporated into skein-theoretic and state-sum invariants.

 Let \(L_+\), \(L_-\), and \(L_0\) denote the usual oriented HOMFLYPT skein triple. For stuck crossings, let \(L_\ast^+\) and \(L_\ast^-\) denote the two oriented stuck crossings, according to the local over/under information.

\begin{figure}[ht]
\centering
\includegraphics[width=0.3\textwidth]{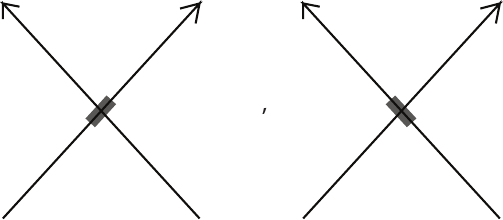}
\caption{The oriented stuck crossings \(L_{\ast}^{+}\) and \(L_{\ast}^{-}\).}
\label{fig:stuck_oriented}
\end{figure}

\begin{definition}[Rigid HOMFLYPT polynomial \cite{Diamantis2026StuckKnots}]
The \emph{rigid HOMFLYPT polynomial} \(P_R\) is defined recursively by the classical HOMFLYPT skein relation
\[
aP_R(L_+) - a^{-1}P_R(L_-) = zP_R(L_0),
\]
together with the stuck crossing relations
\[
P_R(L_\ast^+) = tP_R(L_0)+rP_R(L_+),
\]
\[
P_R(L_\ast^-) = tP_R(L_0)+r^{-1}P_R(L_-),
\]
where \(a,z,t,r\) are independent indeterminates. The normalization is
\[
P_R(\bigcirc)=1,
\qquad
P_R(L\sqcup \bigcirc)
=
\frac{a-a^{-1}}{z}P_R(L).
\]
\end{definition}

The variables \(a\) and \(z\) record the usual HOMFLYPT skein behavior. The new variables record the stuck structure: \(t\) records the smoothing
contribution of a stuck crossing, while \(r\) records the oriented stuck crossing contribution. When no stuck crossings are present, the defining relations reduce to the classical HOMFLYPT polynomial.

\begin{theorem}[\cite{Diamantis2026StuckKnots}]
The polynomial \(P_R\) is invariant under stuck isotopy.
\end{theorem}

This invariant detects the rigidity structure of the diagram. If a stuck knot has \(k\) stuck crossings, then the powers of the additional variables \(t\)
and \(r\) appearing in \(P_R\) are bounded by \(k\). Thus the polynomial carries a natural grading controlled by the number of stuck crossings.

\begin{example}[A stuck curl]
Let \(D\) be a diagram consisting of a single positive stuck crossing \(L_\ast^+\) on a trivial loop. The stuck skein relation gives
\[
P_R(D)=tP_R(L_0)+rP_R(L_+).
\]
Here \(L_0\) is the two-component unlink, while \(L_+\) is a positive classical kink and hence represents the unknot. Therefore
\[
P_R(D)
=
t\left(\frac{a-a^{-1}}{z}\right)+r.
\]
Thus the invariant distinguishes the stuck curl from the classical unknot, even though the underlying classical knot type becomes trivial after forgetting the stuck structure.
\end{example}

\begin{remark}
This construction is related to the tangle insertion philosophy: a stuck crossing is replaced by a controlled linear combination of local classical
configurations. Unlike a pseudo crossing, however, a stuck crossing is not interpreted as ambiguous. The replacement records the persistence of a constrained classical crossing.
\end{remark}


\subsection{A bracket-type state-sum invariant}

There is also an unoriented state-sum construction extending the Kauffman bracket. Classical crossings admit the usual \(A\)- and \(B\)-smoothings. At a
stuck crossing, one introduces a persistent stuck state \(V\), which records the fact that the local rigidity constraint remains visible in the state
expansion.

\begin{figure}[ht]
\centering
\includegraphics[width=0.65\textwidth]{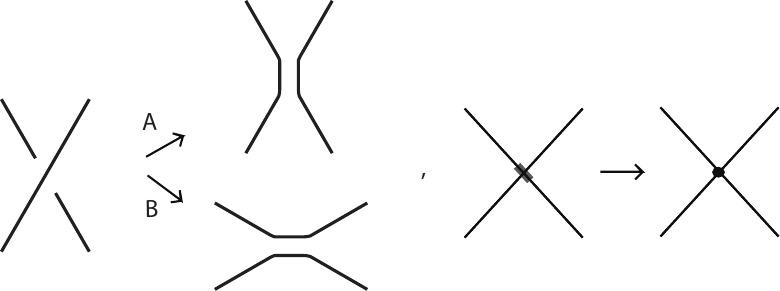}
\caption{Local states at a crossing: \(A\)-smoothing, \(B\)-smoothing, and the
stuck crossing state \(V\).}
\label{fig:states}
\end{figure}

Let \(D\) be a stuck diagram. A state \(s\) is obtained by choosing an \(A\)- or \(B\)-smoothing at each classical crossing and the persistent \(V\)-state at each stuck crossing. Assign the local weights
\[
A\text{-smoothing}: A,\qquad
B\text{-smoothing}: A^{-1},\qquad
V\text{-state}: R,
\]
where \(R\) is an additional indeterminate. If \(\alpha(s)\) and \(\beta(s)\) are the numbers of \(A\)- and \(B\)-smoothings, and if \(\nu(s)\) is the number of stuck states, define
\[
w(s)=A^{\alpha(s)-\beta(s)}R^{\nu(s)}.
\]

\begin{definition}[Stuck bracket \cite{Diamantis2026StuckKnots}]
The \emph{stuck bracket} of a stuck diagram \(D\) is
\[
\langle D\rangle_R
=
\sum_s w(s)\delta^{|s|-1},
\qquad
\delta=-A^2-A^{-2},
\]
where the sum ranges over all states \(s\), and \(|s|\) denotes the number of connected components of the resulting state.
\end{definition}

As in the classical case, the bracket is invariant under the second and third Reidemeister moves, but it must be normalized to account for the first
Reidemeister move. Let \(w(D)\) denote the writhe of \(D\), counting only classical crossings. Stuck crossings do not contribute to the writhe.

\begin{definition}[Normalized stuck bracket \cite{Diamantis2026StuckKnots}]
For a stuck knot \(K^\ast\) represented by a diagram \(D\), define
\[
P_{K^\ast}(A,R)
=
(-A^3)^{-w(D)}\langle D\rangle_R.
\]
\end{definition}

\begin{theorem}
The polynomial \(P_{K^\ast}(A,R)\) is invariant under stuck isotopy.
\end{theorem}

The variable \(R\) records the contribution of persistent stuck states. Hence the state-sum invariant detects rigidity information that disappears after
forgetting the stuck structure.

\begin{example}[Classical curl versus stuck curl]
A positive classical curl has normalized bracket equal to \(1\), as in the classical theory. If the same crossing is declared stuck, then there are no
classical smoothing choices, and one persistent stuck state contributes a factor of \(R\). Since stuck crossings do not contribute to the writhe, the
normalized stuck bracket is
\[
P_{K^\ast}(A,R)=R.
\]
Thus the stuck bracket distinguishes the stuck curl from the classical unknot.
\end{example}


\subsection{Role of stuck knots}

Stuck knots provide a rigidity counterpart to pseudo and singular knot theory (\cite{Bataineh2020,Diamantis2026StuckKnots}). Pseudo knot theory turns a crossing into unresolved information, while singular knot theory turns a crossing into a prescribed double point. Stuck knot theory keeps the crossing classical, but restricts its diagrammatic behavior.

The HOMFLYPT-type and bracket-type invariants show that this rigidity is not merely decorative: it produces algebraic information that is invisible to
classical polynomial invariants after the stuck structure is forgotten. In this way, stuck knot theory fits naturally into the broader landscape of
non-classical knot theories while introducing a distinct principle: local constraints can themselves become topological data.


\section{Bonded Knots and Tied Links}

The theories discussed in the previous sections modify crossings, ambient spaces, or allowed diagrammatic moves. Bonded knots and tied links represent a different direction of generalization (\cite{AicardiJuyumaya2016TiedLinks,DiamantisKauffmanLambropoulou2025Bonded}). In these theories the underlying link diagram is enriched by additional relational data. In the case of bonded links, this extra data is geometric: one adds embedded arcs, called \emph{bonds}, whose endpoints lie on the link. In the case of tied links, the extra data is combinatorial: one records relations among components or strands by means of non-embedded ties.

This distinction is central. A bond is part of the embedded topology of the object and therefore interacts with the diagram through isotopy restrictions, forbidden moves, and spatial-graph phenomena (\cite{Kauffman1989Graphs,
DiamantisKauffmanLambropoulou2025Bonded}). A tie, by contrast, records an equivalence relation among components or strands and is allowed to move flexibly through the diagram. Thus, bonded knots and tied links both enrich classical knot theory by adding relational structure, but they do so in fundamentally different ways (\cite{AicardiJuyumaya2016TiedLinks, Diamantis2021TiedLinksSettings}).


\subsection{Bonded diagrams and types of bonds}

Informally, a bonded link is a link together with a finite collection of auxiliary arcs whose endpoints lie on the link. These arcs are the bonds.

\begin{definition}[Bonded link]
A \emph{bonded link} is a pair \((L,B)\), where \(L\) is a link in \(S^3\) and \(B\) is a finite collection of disjoint arcs properly embedded in \(S^3\setminus L\), such that each arc in \(B\) has its two endpoints attached to \(L\). These endpoints are called \emph{nodes}. If \(L\) has one component, then \((L,B)\) is called a \emph{bonded knot}.
\end{definition}

A bonded link diagram is obtained by projecting both the link and the bonds to the plane, with the usual over/under information at double points. Crossings may occur between two link arcs, between a link arc and a bond, or between two bonds. The nodes are the points at which the bonds attach to the link.

\begin{remark}
Bonded links may also be viewed as special embedded trivalent graphs \cite{Kauffman1989Graphs,
DiamantisKauffmanLambropoulou2025Bonded}. The ordinary edges are the link arcs, while the bond edges are distinguished auxiliary arcs. This graph-theoretic viewpoint is useful, but it is important that bonds are not simply ordinary edges of the original knot or link. They represent additional structure attached to the link.
\end{remark}

There are several levels of restriction on the allowed bonds.

\begin{definition}[Long, standard, and tight bonds \cite{DiamantisKauffmanLambropoulou2025Bonded}]
A \emph{long bond} is an embedded bond arc that may itself be knotted or linked with the rest of the diagram.

A \emph{standard bond} is an unknotted and unlinked bond, usually presented as a simple segment connecting two nearby strands in an \(H\)-neighborhood.

A \emph{tight bond} is a standard bond that does not interact with link arcs; equivalently, the bond lies in an \(H\)-neighborhood that contains no additional crossings involving the bond.
\end{definition}

\begin{figure}[ht]
    \centering
    \includegraphics[width=0.75\textwidth]{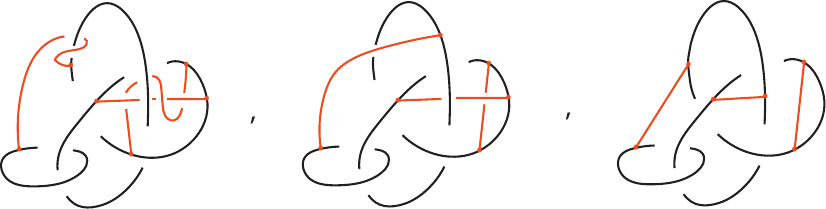}
    \caption{Schematic forms of long, standard and tight bonds.}
    \label{fig:standard-tight-bonds}
\end{figure}

Thus, long bonded links form the most general category. Standard bonded links are obtained when the bonds themselves have been simplified, while tight bonded links represent the most local version of the theory. For an illustration see Figure~\ref{fig:standard-tight-bonds}. These distinctions are useful because different categories lead to different move systems and different invariant constructions.

A useful feature of the theory is that, in the diagrammatic settings considered in \cite{DiamantisKauffmanLambropoulou2025Bonded}, long bonds can be related to standard forms, and standard bonds can be brought into tight form by removing interactions between the bond and the link arcs. Thus, the three categories are not unrelated; rather, they form a hierarchy of diagrammatic models for the same general idea of a link equipped with embedded auxiliary connections.


\subsection{Bonded isotopy}

The equivalence relation for bonded links must preserve both the underlying link and the bond structure. There are two natural versions of this equivalence.

\begin{definition}[Topological and rigid bonded isotopy]
Two bonded links \((L_1,B_1)\) and \((L_2,B_2)\) are equivalent by \emph{topological vertex isotopy} if there is an ambient isotopy of \(S^3\) taking \(L_1\) to \(L_2\) and \(B_1\) to \(B_2\), preserving the attachment of bonds to link arcs.

They are equivalent by \emph{rigid vertex isotopy} if, in addition, the local neighborhoods of the nodes are preserved as rigid vertex neighborhoods.
\end{definition}

\begin{remark}
Rigid vertex isotopy is stricter than topological vertex isotopy. This distinction is important because some invariants are natural in the topological category, while others, such as tangle insertion invariants, are well defined only in the rigid category.
\end{remark}

On the level of diagrams, bonded isotopy is described by Reidemeister-type moves for the underlying link together with additional moves involving bonds and their nodes. In the full theory of long and standard bonded links, these moves include Reidemeister moves on bond arcs, mixed moves between bond arcs and link arcs, vertex slide moves, and vertex twist moves. The precise list depends on whether one works in the topological vertex category or in the rigid vertex category.

For the purposes of this survey, we shall not list the complete move systems for all categories (the interested reader is referred to \cite{DiamantisKauffmanLambropoulou2025Bonded}). Instead, we emphasize the hierarchy
\[
\text{long bonded links}
\longrightarrow
\text{standard bonded links}
\longrightarrow
\text{tight bonded links}.
\]
The tight category is the most local and is the one most closely connected with the bonded bracket polynomial and the braid-theoretic formulation below. In this category, each bond is contained in an \(H\)-neighborhood and does not cross any link arc. Thus, the bond behaves like a local rigid auxiliary connection rather than a long embedded arc moving through the diagram.

\begin{figure}[ht]
    \centering
    \includegraphics[width=0.65\textwidth]{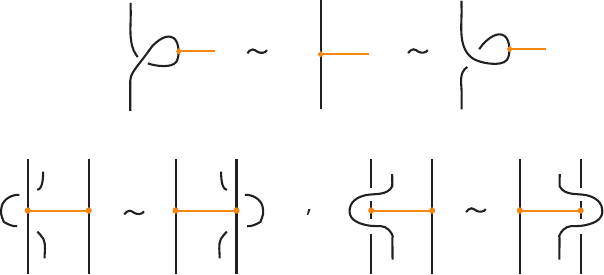}
    \caption{Typical local behavior of a tight bond. In the tight category the bond is kept inside an \(H\)-neighborhood and does not cross link arcs.}
    \label{fig:tight-bond-moves}
\end{figure}

\begin{figure}[ht]
    \centering
    \includegraphics[width=0.75\textwidth]{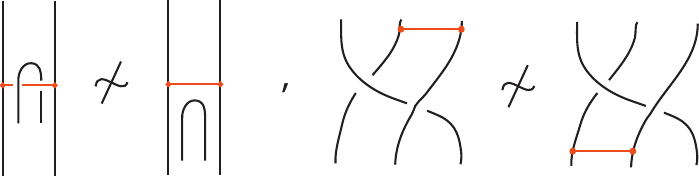}
    \caption{Some forbidden moves in the theory of bonded links.}
    \label{fig:forb-moves}
\end{figure}

Since a bond is an embedded arc, it cannot pass freely
through strands, through other bonds, or through crossings unless the motion is realized by the allowed bonded moves (see Figure~\ref{fig:forb-moves}). In particular, a bond cannot be treated as a phantom relation that may be pulled arbitrarily through the diagram. This embedded nature is exactly what gives bonded knot theory its spatial-graph character.

\subsection{Invariants of bonded links}

The construction of invariants for bonded links must account for two layers of information: the isotopy type of the underlying link and the placement of the bonding structure. Several techniques are available, depending on the category of bonded links under consideration.

For topological bonded links, one useful method is \emph{unplugging}. Since a bonded link may be regarded as an embedded graph, one can remove selected edges
at the trivalent nodes and obtain a collection of ordinary classical links. For rigid bonded links, a different method is more natural: \emph{tangle insertion}. One replaces each bond by a small band and inserts a chosen \(2\)-tangle into that band. The resulting object is an ordinary classical link, to which classical invariants may be applied.

\begin{theorem}[Tangle insertion invariant]
Let \((L,B)\) be a standard bonded link in the rigid vertex category. Replacing each bond by a band and inserting a fixed family of \(2\)-tangles produces classical links whose isotopy classes depend only on the rigid vertex isotopy class of \((L,B)\). Therefore, applying any classical link invariant after tangle insertion yields an invariant of rigid bonded links.
\end{theorem}

This method is specific to the rigid setting. In the topological category, the nodes may rotate in ways that make a fixed tangle insertion ill defined under isotopy.

A more explicit skein-theoretic construction is given by the bonded bracket polynomial. For rigid tight bonded links, one assigns local skein relations not only at ordinary crossings, as in the Kauffman bracket, but also at bonds.

\begin{definition}[Bonded bracket polynomial]
Let \(L\) be a rigid tight bonded link diagram. The \emph{bonded bracket polynomial} is a Kauffman bracket type polynomial obtained by applying the usual Kauffman bracket relations at classical crossings, together with the following additional local replacement rule at each bond:
\[
L_b = aL_0 + bL_+ + bL_-.
\]
Here \(L_b\) denotes the local diagram consisting of an \(H\)-neighborhood of a bond, \(L_+\) and \(L_-\) denote the corresponding positive and negative classical crossings obtained by replacing the bond locally, and \(L_0\) denotes the corresponding smoothing. The diagrams \(L_b,L_0,L_+\), and \(L_-\) are identical outside the chosen \(H\)-neighborhood.
\end{definition}

\begin{theorem}
The bonded bracket polynomial is invariant under regular rigid vertex isotopy in the category of tight bonded links.
\end{theorem}

The bonded bracket shows that bonds are not merely decorative additions to a diagram. They contribute algebraic information that may distinguish bonded knots with the same underlying classical knot type. In this sense, bonded invariants detect both the topology of the backbone and the topology of the additional embedded interactions.


\subsection{Bonded braids, enhanced bonds, and bonded knotoids}

\begin{figure}[ht]
    \centering
    \includegraphics[width=0.55\textwidth]{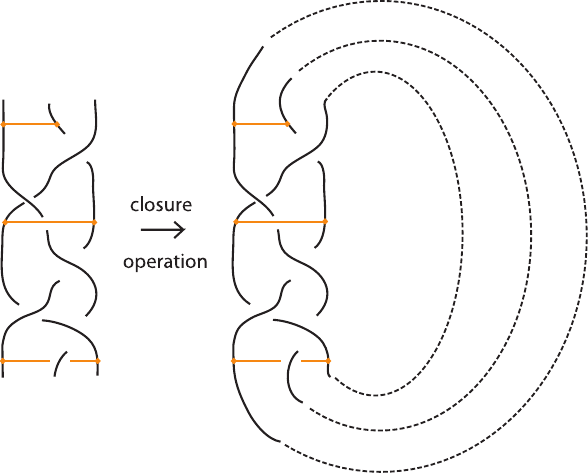}
    \caption{A bonded braid and its closure.}
    \label{fig:bonded-braid}
\end{figure}

As in classical knot theory, the diagrammatic theory of bonded links has a braid-theoretic counterpart. A bonded braid is a braid equipped with bonds connecting points on its strands as illustrated in Figure~\ref{fig:bonded-braid}. The closure of a bonded braid is obtained by closing the braid in the usual way while retaining the bonds.

\begin{definition}[Bonded braid \cite{DiamantisKauffmanLambropoulou2025Bonded}]
A \emph{bonded braid on \(n\) strands} is a classical braid on \(n\) strands together with a finite collection of disjoint embedded horizontal arcs, called
bonds, whose endpoints lie on the braid strands. A bond joining the \(i\)-th and \(j\)-th strands is denoted by \(b_{i,j}\). A bond joining two consecutive
strands \(i\) and \(i+1\) is called an elementary bond and is denoted by \(b_i\).
\end{definition}

The algebraic structure associated with tight bonded braids is the bonded braid monoid \(BB_n\). It is generated by the ordinary braid generators \(\sigma_1^{\pm1},\ldots,\sigma_{n-1}^{\pm1}\) together with elementary bond generators
\[
b_1,\ldots,b_{n-1}.
\]
The generators \(\sigma_i\) satisfy the usual braid relations and the elementary bond generators satisfy the following relations with each other and with the braid generators:
\[
b_i b_j=b_j b_i,
\quad \text{and} \quad
b_i\sigma_j=\sigma_j b_i
\qquad \text{for}\ |i-j|>1,
\]
\[
b_i\sigma_i=\sigma_i b_i,
\qquad 1\leq i\leq n-1,
\]
\[
b_i\sigma_{i+1}\sigma_i
=
\sigma_{i+1}\sigma_i b_{i+1}
\quad \text{and} \quad
\sigma_i\sigma_{i+1}b_i
=
b_{i+1}\sigma_i\sigma_{i+1},
\qquad \text{for}\ 1\leq i\leq n-2.
\]

A central feature of this algebraic structure is its relation with the singular braid monoid.

\begin{theorem}[\cite{DiamantisKauffmanLambropoulou2025Bonded}]
The bonded braid monoid \(BB_n\) is isomorphic to the singular braid monoid. Under this correspondence, the elementary bond generator \(b_i\) plays the formal algebraic role of the singular braid generator \(\tau_i\).
\end{theorem}

This result should be interpreted algebraically, not geometrically. A bond is not a singular crossing; it is an embedded auxiliary arc connecting two strands. Nevertheless, there is a useful local degeneration picture behind the
isomorphism. In the tight category, an elementary bond \(b_i\) is confined to a small \(H\)-neighborhood connecting two adjacent strands. If this neighborhood is contracted, the bond is visually collapsed to a small contact point between the strands. The resulting limiting picture is not a transverse singular crossing, but rather a tangential contact-type singularity. At the level of local braid relations, however, this contracted bond satisfies the same formal
relations as the singular braid generator \(\tau_i\). Thus the correspondence \(b_i\leftrightarrow \tau_i\) records a formal algebraic analogy, while the geometric meanings remain distinct.

Bonded braids also satisfy analogues of the Alexander and Markov theorems.

\begin{theorem}[Alexander theorem for bonded links \cite{DiamantisKauffmanLambropoulou2025Bonded}]
Every oriented topological standard bonded link can be represented, up to topological vertex isotopy, as the closure of a standard or tight bonded braid.
\end{theorem}

The Markov-type theorem for bonded braids is most naturally formulated in terms of bonded \(L\)-moves. These moves play the role of stabilization and conjugation in the classical Markov theorem, but must be supplemented by bond-commuting moves, since the position of embedded bonds is not fully controlled by ordinary \(L\)-moves.

\begin{theorem}[Bonded \(L\)-move Markov theorem \cite{DiamantisKauffmanLambropoulou2025Bonded}]
Two bonded braids have topologically equivalent closures if and only if they are related by bonded braid isotopy, bonded \(L\)-moves, and elementary
bond-commuting moves
\[
\alpha b_i \sim b_i\alpha,
\qquad
\alpha\in BB_n.
\]
\end{theorem}

Equivalently, since \(L\)-moves encode the geometric content of stabilization and conjugation, this theorem may be regarded as the \(L\)-move formulation of the Markov theorem for bonded braids. A fully algebraic Markov-style presentation can be obtained by replacing the \(L\)-moves with the corresponding stabilization and conjugation moves, together with the necessary bond-commuting moves.

Bonded knot theory also admits enhanced versions. One may introduce two types of bonds, often interpreted as attracting and repelling interactions. Algebraically, these two bond types may be treated as inverses. This turns the bonded braid monoid into an enhanced bonded braid group.

\begin{definition}[Enhanced bonded braid group \cite{DiamantisKauffmanLambropoulou2025Bonded}]
The \emph{enhanced bonded braid group} is obtained from the bonded braid monoid by adjoining inverses to the bond generators, corresponding to two mutually inverse types of enhanced bonds.
\end{definition}

This extension is useful when bonds are meant to model interactions of different signs or types, such as attractive and repulsive constraints.


\subsection{Tied links and tied braids}

Tied links form a related but distinct theory \cite{AicardiJuyumaya2016TiedLinks}. A tie is not an embedded arc in the complement of a link. Rather, it records a relation among components, and the set of ties determines a partition of the set of components.

\begin{definition}[Tied link]
An \emph{oriented tied link} \(L(P)\) consists of an oriented link \(L\) in \(S^3\) together with a set \(P\) of ties. Each tie connects two points on the
components of \(L\), and the collection of ties determines a partition of the set of components of \(L\). If \(P=\emptyset\), then \(L(P)\) is an ordinary
oriented link.
\end{definition}

The important point is that ties are not embedded arcs. They may be drawn as springs or line segments, but they are allowed to pass through the diagram.
Their role is to encode the equivalence relation generated by the statement that two components are tied.

\begin{definition}[Tie isotopy]
Two tied links \(L(P)\) and \(L'(P')\) are \emph{tie isotopic} if the underlying links \(L\) and \(L'\) are ambient isotopic and the sets of ties \(P\) and \(P'\) define the same partition of the set of components.
\end{definition}

Thus, one may introduce or delete redundant ties without changing the tied link, provided the induced partition of components is unchanged. For example,
an additional tie between two components already connected through a chain of ties does not change the tied link.

The braid-theoretic counterpart is the tied braid monoid.

\begin{definition}[Tied braid monoid \cite{AicardiJuyumaya2000BraidsTies,AicardiJuyumaya2016TiedLinks}]
The \emph{tied braid monoid} \(TM_n\) is generated by the braid generators \(\sigma_1^{\pm1},\ldots,\sigma_{n-1}^{\pm1}\) and tie generators \(\eta_1,\ldots,\eta_{n-1}\). The generators \(\sigma_i\) satisfy the usual braid relations. The tie generators satisfy idempotent and compatibility relations, including
\[
\eta_i^2=\eta_i,
\qquad
\eta_i\eta_j=\eta_j\eta_i,
\]
\[
\eta_i\sigma_i=\sigma_i\eta_i,
\qquad
\eta_i\sigma_j=\sigma_j\eta_i
\quad (|i-j|>1),
\]
together with the standard adjacent braid--tie relations that allow ties to move through braid generators.
\end{definition}

Generalized ties \(\eta_{i,j}\), joining the \(i\)-th and \(j\)-th strands, are obtained by conjugating elementary ties by braid generators. They are
transparent with respect to the strands between their endpoints and possess an elasticity property. These features imply the fundamental mobility property of
tied braids.

\begin{proposition}[Mobility property \cite{AicardiJuyumaya2016TiedLinks}]
Every tied braid can be written, up to the defining relations of \(TM_n\), as a product
\[
\alpha=\beta\gamma,
\]
where \(\beta\in B_n\) is an ordinary braid and \(\gamma\) is a product of generalized ties. Moreover, the product of generalized ties determines a partition of the set of strands.
\end{proposition}

This mobility property is the key feature distinguishing ties from bonds. A bond cannot in general be moved freely through the diagram, because it is an embedded arc. A tie can be moved to the top or bottom of a braid because it only records partition data.

The Alexander and Markov theorems extend to tied links.

\begin{theorem}[Alexander theorem for tied links \cite{AicardiJuyumaya2016TiedLinks}]
Every oriented tied link is tie isotopic to the closure of a tied braid.
\end{theorem}

\begin{theorem}[Markov theorem for tied braids \cite{AicardiJuyumaya2016TiedLinks}]
Two tied braids have tie isotopic closures if and only if they are related by a finite sequence of braid relations and the following moves:
\[
\begin{array}{rcll}
\text{Conjugation:}
&
\alpha\beta
&\sim&
\beta\alpha,
\qquad
\alpha,\beta\in TM_n,
\\[0.3cm]
\text{Stabilization:}
&
\alpha
&\sim&
\alpha\sigma_n^{\pm1},
\qquad
\alpha\in TM_n,
\\[0.3cm]
\text{Tie stabilization:}
&
\alpha
&\sim&
\alpha\eta_{i,j},
\qquad
s_\alpha(i)=j.
\end{array}
\]
Here \(s_\alpha\) denotes the permutation associated to the braid obtained from \(\alpha\) by forgetting its ties.
\end{theorem}

From an algebraic perspective, tied links lead naturally to braid-and-tie algebras, or \(bt\)-algebras. These algebras extend braid-related algebras by adjoining idempotent tie generators and compatibility relations with the braid generators. They support trace constructions and therefore provide a natural setting for constructing polynomial invariants of tied links \cite{AicardiJuyumaya2000BraidsTies,AicardiJuyumaya2016TiedLinks,
AicardiJuyumaya2018KauffmanTied}.


\subsection{Tied pseudo and tied singular extensions}

Tied links can be combined with the crossing-based theories discussed earlier. A \emph{tied singular link} is a singular link equipped with ties, while a \emph{tied pseudo link} is a pseudo link equipped with ties \cite{AicardiJuyumaya2018TiedSingular,Diamantis2021TiedPseudoLinks}. Thus, tied singular links combine prescribed transverse double points with partition data, whereas tied pseudo links combine unresolved crossings with partition data.

The braid-theoretic structures reflect this combination.

\begin{definition}[Tied pseudo braid monoid \cite{Diamantis2021TiedPseudoLinks}]
The \emph{tied pseudo braid monoid} \(TPM_n\) is generated by
\[
\sigma_1^{\pm1},\ldots,\sigma_{n-1}^{\pm1},
\qquad
p_1,\ldots,p_{n-1},
\qquad
\eta_1,\ldots,\eta_{n-1}.
\]
The generators satisfy the defining relations of the pseudo braid monoid \(PM_n\), the defining relations of the tied braid monoid \(TM_n\), and additional compatibility relations describing how ties interact with pre-crossings.
\end{definition}

Similarly, the tied singular braid monoid \(TSM_n\) is generated by the classical braid generators, singular generators \(\tau_i\), and tie generators \(\eta_i\), with analogous compatibility relations.

As in the untied case, the pseudo and singular versions have the same formal monoid structure after replacing singular crossings by pre-crossings.

\begin{theorem}
There is an isomorphism
\[
TSM_n\longrightarrow TPM_n
\]
defined on generators by
\[
\sigma_i^{\pm1}\longmapsto\sigma_i^{\pm1},
\qquad
\tau_i\longmapsto p_i,
\qquad
\eta_i\longmapsto\eta_i.
\]
\end{theorem}

This isomorphism is formal and algebraic. It does not identify the geometric meaning of a singular crossing with that of a pre-crossing. Rather, it says
that the monoid relations are compatible with the replacement \(\tau_i\leftrightarrow p_i\), while the tie generators remain unchanged.

The Alexander theorem also extends to tied pseudo links.

\begin{theorem}[Alexander theorem for tied pseudo links \cite{Diamantis2021TiedPseudoLinks}]
Every oriented tied pseudo link is equivalent to the closure of a tied pseudo braid.
\end{theorem}

The Markov theorem combines the pseudo Markov moves with the tie-stabilization move.

\begin{theorem}[Markov theorem for tied pseudo braids \cite{Diamantis2021TiedPseudoLinks}]
Two tied pseudo braids have equivalent closures if and only if they are related by a finite sequence of tied pseudo braid relations and the following moves:
\[
\begin{array}{rcll}
\text{Commuting:}
&
\alpha p_i
&\sim&
p_i\alpha,
\qquad
\alpha\in TPM_n,
\\[0.3cm]
\text{Conjugation:}
&
\beta
&\sim&
\gamma^{\pm1}\beta\gamma^{\mp1},
\qquad
\beta\in TPM_n,\ \gamma\in B_n,
\\[0.3cm]
\text{Real stabilization:}
&
\alpha
&\sim&
\alpha\sigma_n^{\pm1},
\qquad
\alpha\in TPM_n,
\\[0.3cm]
\text{Pseudo-stabilization:}
&
\alpha
&\sim&
\alpha p_n,
\qquad
\alpha\in TPM_n,
\\[0.3cm]
\text{Tie stabilization:}
&
\alpha
&\sim&
\alpha\eta_{i,j},
\qquad
s_\alpha(i)=j.
\end{array}
\]
Here \(s_\alpha\) denotes the permutation associated to the braid obtained from \(\alpha\) by forgetting its ties and pre-crossings when appropriate.
\end{theorem}

Thus, tied pseudo links combine two independent kinds of non-classical data: unresolved local crossing information and global partition data among components or strands. The former is represented algebraically by the generators \(p_i\), while the latter is represented by the idempotent tie generators \(\eta_i\).


\subsection{Bonds versus ties}

Bonded knots and tied links are naturally compared because both enrich classical diagrams with additional relational structure. However, the comparison should not obscure the fundamental difference between them.

A bond is an embedded geometric arc. It has a position in space, may be knotted or linked with the diagram, may interact with link arcs, and is constrained by the allowed isotopy moves. Because of this, bonded links are closely related to spatial graph theory and require careful treatment of topological versus rigid vertex isotopy, forbidden moves, and embedded-arc invariants.

A tie is not an embedded arc. It is a diagrammatic notation for an equivalence relation among components or strands. Ties are transparent and elastic; they may pass through the diagram and can be moved to the top or bottom of a tied braid. Algebraically, this leads to idempotent tie generators, generalized ties, mobility properties, and braid-and-tie algebras.

\begin{center}
\tiny
\[
\begin{array}{c|c|c}
 & \text{Bonded links} & \text{Tied links} \\
\hline
\text{extra data} & \text{embedded arcs} & \text{partition data} \\
\text{local behavior} & \text{geometrically constrained} & \text{transparent and mobile} \\
\text{main algebraic generators} & b_i & \eta_i \\
\text{basic algebraic feature} & BB_n\cong SM_n & \eta_i^2=\eta_i \\
\text{geometric nature} & \text{spatial-graph-like} & \text{braid-and-tie algebraic} \\
\text{typical invariants} & \text{unplugging, tangle insertion, bonded bracket}
& \text{trace invariants from } bt\text{-algebras}
\end{array}
\]
\end{center}

This comparison places bonded knots and tied links within the broader landscape of non-classical knot theory. Pseudo knots modify crossing information through ambiguity; singular knots introduce prescribed double points; stuck knots impose rigidity constraints; virtual knots modify the diagrammatic or ambient setting. Bonded knots and tied links instead enrich the diagram with auxiliary relational data external to the crossings themselves.

This shift is conceptually important. It shows that non-classical knot theory is not only about changing crossings, move systems, or ambient spaces. It also includes theories in which knots and links carry additional structures recording interactions, constraints, or relations among different parts of the diagram. Bonded knots emphasize the topology of embedded interactions, while tied links emphasize the algebra of partition data. Together they provide two complementary models for relational structure in knot theory.


\section{Further Diagrammatic and Ambient Generalizations}

The theories discussed so far modify classical knot theory by changing local crossing data, imposing rigidity constraints, altering the ambient three-manifold, or enriching diagrams with additional relational structure. We now briefly recall several further generalizations in which the classical diagrammatic framework is modified in different ways. Virtual and welded knot theories arise from relaxing planarity and changing the allowed move system, while knotoids and braidoids arise by replacing closed diagrams with open-ended ones (\cite{Kauffman1999,FennRimanyiRourke1997BraidPermutation,Turaev2012,
GugumcuLambropoulou2021Braidoids}).


\subsection{Virtual knots}

Virtual knot theory, introduced by Kauffman~\cite{Kauffman1999}, extends classical knot theory by allowing diagrams to contain, in addition to classical crossings, \emph{virtual crossings}. A virtual crossing is usually represented by a circled intersection and carries no over/under information.

It is important, however, not to interpret virtual crossings as pseudo crossings. A pseudo crossing represents missing classical crossing information; a virtual crossing records the failure of a planar drawing to represent the underlying object without auxiliary intersections. In this sense, virtual knot theory changes the diagrammatic and ambient setting rather than simply introducing another kind of unresolved crossing.

\begin{definition}[Virtual knot diagram \cite{Kauffman1999}]
A \emph{virtual knot diagram} is a knot diagram with two types of crossings: classical crossings, which carry over/under information, and virtual crossings, usually drawn as circled double points. A \emph{virtual knot} is an equivalence class of virtual knot diagrams under the classical Reidemeister moves together with the virtual Reidemeister moves encoded by the detour principle. The forbidden moves, in which a strand is moved across a virtual crossing as if it were a classical crossing, are not allowed.
\end{definition}

The virtual Reidemeister moves are often summarized by the \emph{detour principle}: any arc that meets the rest of the diagram only in virtual crossings may be replaced by another such arc, with all new intersections declared virtual. This principle expresses the fact that virtual crossings are not genuine local crossings of the knot, but bookkeeping devices introduced by the planar representation.

The geometric interpretation of virtual knots is given by knots embedded in thickened surfaces \(\Sigma\times[0,1]\), considered up to isotopy and stabilization of the surface. This interpretation, developed by Carter, Kamada, and Saito~\cite{CarterKamadaSaito2002}, shows that virtual knot theory may be viewed as the study of knots in thickened surfaces by means of planar diagrams. Kuperberg's theorem~\cite{Kuperberg2003} further states that each virtual knot has a minimal-genus surface representative, unique up to homeomorphism.

Virtual knots also admit a natural Gauss diagram formulation. In the classical case, not every Gauss diagram is realizable by a planar knot diagram. In the virtual setting, every Gauss diagram is realizable. Thus, virtual knot theory can also be regarded as a natural completion of classical knot theory from the viewpoint of Gauss diagrams.

The braid-theoretic counterpart is given by the virtual braid groups \(VB_n\), generated by the classical braid generators \(\sigma_i\) together with virtual generators \(v_i\) \cite{KauffmanLambropoulou2006VirtualBraids}. Virtual braids admit closure operations and satisfy Alexander and Markov type theorems. However, no single Hecke-type quotient
plays the same canonical role for virtual knots as the Hecke algebra does in classical knot theory. Consequently, virtual knot invariants are often constructed using Gauss diagrams, parity methods, quandles, biquandles, state-sum models, and refinements such as the arrow polynomial.

Virtual knot theory therefore occupies a central position among non-classical knot theories. It preserves a Reidemeister-type diagrammatic calculus while changing the surface on which the diagram is naturally represented.


\subsection{Welded knots}

Welded knot theory is obtained from virtual knot theory by allowing one of the two forbidden moves. It is therefore a quotient theory of virtual knot theory: some virtual information is lost, but the resulting theory has important algebraic and geometric interpretations \cite{Kamada2007WeldedBraids}.

\begin{definition}[Welded knot]
A \emph{welded knot} is an equivalence class of virtual knot diagrams under the classical and virtual Reidemeister moves together with one of the forbidden moves, usually called the over-forbidden move.
\end{definition}

The asymmetry between the two forbidden moves is significant. Allowing both forbidden moves collapses the theory much further, while allowing only the over-forbidden move produces the welded category. This category retains enough structure to support a rich braid theory, but is flexible enough to connect with higher-dimensional topology.

The braid-theoretic counterpart of welded knots is given by welded braid groups, which are closely related to loop braid groups and groups of motions of circles in three-space \cite{FennRimanyiRourke1997BraidPermutation,Damiani2017LoopBraidGroups}. Geometrically, welded knots are related to ribbon torus knots in four-dimensional topology through Satoh's Tube construction \cite{Satoh2000}. Under this construction, a welded knot diagram gives rise to a knotted torus in \(S^4\), providing a bridge between virtual diagrammatics and higher-dimensional knot theory.

From the perspective of this survey, welded knots illustrate how a controlled change in the allowed move system can produce a new and meaningful quotient theory. They sit between virtual knot theory, braid group theory, and four-dimensional diagrammatic topology.


\subsection{Knotoids and braidoids}

Another direction of generalization changes the global form of the diagram rather than the local crossing type. Classical knot diagrams are closed. Knotoid diagrams, introduced by Turaev~\cite{Turaev2012}, are open-ended diagrams with two distinguished endpoints.

\begin{definition}[Knotoid diagram]
A \emph{knotoid diagram} is an immersed oriented interval in a surface, with classical crossing information at double points, considered up to Reidemeister moves performed away from the endpoints.
\end{definition}

The endpoints are part of the structure. They cannot simply be pulled across strands, and this restriction gives rise to the characteristic forbidden endpoint moves of knotoid theory. As a result, knotoids can detect information about open curves that would disappear after an arbitrary closure.

Knotoids are closely related to classical knots, but they are not merely ``open knots''. Different closure operations, such as over-closure and under-closure, may produce classical knots or links from a knotoid, but the knotoid itself retains information about the position and interaction of the endpoints before closure. This makes knotoids particularly useful in settings where one wants to study open chains without imposing an artificial closure.

Braidoids provide the corresponding braid-like objects \cite{GugumcuLambropoulou2021Braidoids}. A braidoid diagram is a braid-type diagram that may contain strands with free endpoints. Braidoids play for knotoids a role analogous to the role played by braids for classical links: they provide an algebraic and diagrammatic framework in which open-ended knot-like objects can be represented.

As in the classical case, one has Alexander and Markov type results relating knotoids and braidoids. The Alexander-type theorem states that knotoids can be represented by braidoids, while the corresponding Markov-type theorem describes when two braidoids determine equivalent knotoids. These results show that the braid-theoretic philosophy extends beyond closed diagrams, although the presence of endpoints requires additional care.

The knotoid and braidoid framework also admits pseudo and ambient generalizations. Pseudo knotoids and pseudo braidoids were introduced in \cite{Diamantis2021TiedPseudoLinks}, extending knotoid theory by allowing unresolved crossings together with endpoint-sensitive equivalence. These
objects combine the local ambiguity of pseudo knot theory with the global endpoint structure of knotoids. The theory was further developed in annular and toroidal settings, where the presence of a nontrivial ambient surface or solid-torus structure introduces additional winding information
\cite{Diamantis2022KnotoidsTorus,
DiamantisLambropoulouMahmoudi2024ToroidalPseudoKnots,
DiamantisLambropoulouMahmoudi2025ToroidalKnotoids}.

More recently, these ideas were extended to doubly periodic diagrammatic settings. In particular, doubly periodic pseudo tangles provide a periodic analogue of pseudo diagrammatic theory, while the equivalence theory of doubly periodic tangles contains, as special cases or limiting diagrammatic settings, many of the knot-theoretic structures discussed in this survey
\cite{DiamantisLambropoulouMahmoudi2024DoublyPeriodicPseudoTangles,
DiamantisLambropoulouMahmoudi2026DoublyPeriodicTangles}.

Knotoids and braidoids therefore represent a different kind of non-classical generalization. They do not primarily modify crossings, add bonds, or change the ambient surface. Instead, they modify the global diagrammatic object by allowing open-ended curves with endpoint-sensitive equivalence.


\subsection{Other directions}

Many further variants of knot theory fit into the same general philosophy. Flat virtual knots are obtained from virtual knots by forgetting the over/under information at classical crossings, while free knots go further by considering diagrammatic structures with even less crossing information \cite{Kauffman1999,Manturov2009FreeKnots}. Spatial graph theory replaces embedded circles by embedded graphs, while singular and bonded graph theories enrich such graphs with additional local or relational data.

Other extensions arise by changing the equivalence relation, adding decorations to components or crossings, considering links in other ambient manifolds, or allowing diagrams on surfaces with additional structure. Each such theory
modifies one or more ingredients of the classical framework: the crossing data, the allowed local moves, the ambient topology, the closure condition, or the auxiliary data carried by the diagram.

The common feature is that classical knot theory serves as a structural template. Once one identifies which part of the classical framework is being modified, the natural questions reappear: what are the appropriate diagrams, what are the allowed moves, is there a braid-theoretic model, and which invariants survive or must be replaced?


\section{Concluding Remarks and Open Problems}

The aim of this survey has been to describe several ways in which classical knot theory extends beyond its original setting. Classical knot theory provides a stable structural template built from diagrams, Reidemeister moves, braid representations, skein relations, and algebraic invariants. The non-classical theories discussed here arise by modifying one or more of these ingredients.

The examples considered in this survey illustrate different types of modification. Skein modules extend polynomial invariants from links in \(S^3\) to links in more general three-manifolds. Pseudo and singular knot theories modify the local crossing structure: pseudo crossings encode unresolved classical crossing information, while singular crossings encode prescribed transverse double points. Stuck knot theory keeps the crossing classical, but adds a rigidity constraint that restricts its behavior under isotopy. Bonded knots and tied links shift the focus away from crossings by enriching diagrams with auxiliary relational data. Virtual and welded knots arise from changing the diagrammatic and ambient  setting, while knotoids and braidoids replace closed diagrams by open-ended ones.

Taken together, these theories show that non-classical knot theory is not a collection of unrelated extensions. Rather, it forms a broad landscape organized around a small number of recurring structural themes: local crossing data, allowed diagrammatic moves, braid-like representatives, skein-type relations, trace constructions, ambient topology, and auxiliary diagrammatic structure. One of the striking features of this landscape is that much of the classical machinery survives, but often in modified form. Braid groups are replaced by generalized braid monoids or groups, Reidemeister moves are enlarged or restricted, skein relations acquire new local terms, and polynomial invariants
gain additional variables recording ambiguity, singularity, rigidity, bonding, or ambient information.

This perspective suggests several directions for further development.

\begin{itemize}
    \item \textbf{Generalized braid--algebra--trace frameworks.}
    Classical link invariants arise naturally from braid groups, Hecke algebras, and Markov traces. A major open direction is to develop analogous trace constructions for non-classical braid monoids and Hecke-type algebras, especially in settings involving pseudo, singular, bonded, or constrained crossings.

    \item \textbf{Hybrid non-classical theories.}
    Many of the theories discussed here can be combined. One may consider, for example, pseudo-bonded knots, singular-bonded knots, stuck-virtual knots, or bonded knotoids. Such hybrid theories raise natural questions about their move systems, braid representatives, skein relations, and polynomial invariants.

    \item \textbf{Relations among resolution, degeneration, and insertion.}
    Pseudo crossings lead to resolution sets and weighted resolution sets, singular crossings lead to signed skein expansions and finite-type invariants, and bonded or rigid structures lead naturally to tangle insertion. Understanding the precise relationships among these mechanisms may clarify how different non-classical theories encode local replacement data.

    \item \textbf{Skein modules in generalized settings.}
    Skein modules provide a powerful way to study links in three-manifolds. It is natural to ask how skein-module methods extend to diagrams carrying non-classical data, such as pseudo crossings, singular crossings, stuck crossings, bonds, ties, or virtual structure.

    \item \textbf{Probabilistic and statistical invariants.}
    The weighted resolution set of a pseudo knot suggests a probabilistic viewpoint on unresolved crossing data. Similar ideas may be useful for other generalized theories in which uncertainty, ambiguity, or multiple local states are intrinsic to the model.

    \item \textbf{Braid theory beyond closed classical links.}
    Braidoids, bonded braids, virtual braids, welded braids, and mixed braids show that the Alexander--Markov philosophy extends far beyond ordinary
    links in \(S^3\). A systematic comparison of these braid theories may lead to a clearer understanding of which features of classical braid theory are
    structural and which are specific to classical knots.

    \item \textbf{Applications and interpretation.}
    Many non-classical theories are motivated by situations in which classical knots are too restrictive: open chains, molecular bonds, periodic structures, uncertainty in crossings, rigidity constraints, or higher-dimensional interactions. Clarifying the mathematical role of these
    applications can lead both to better models and to new purely topological questions.
\end{itemize}

These problems indicate that the algebraic and geometric aspects of non-classical knots remain far from exhausted. Classical knot theory continues to provide the guiding template, but each generalization reveals phenomena that are invisible in the classical setting. The recurring challenge is to determine which parts of the classical theory survive, which must be modified, and which new structures appear only after leaving the classical framework. In this sense, non-classical knot theory is not merely an extension of classical knot theory, but a systematic exploration of the many ways in which knot-theoretic objects can carry topology, algebra, geometry, and additional structure.

\bibliographystyle{amsplain}
\bibliography{references}

@article{Artin1925, author = {E. Artin}, title = {Theorie der Z{\"o}pfe}, journal = {Abhandlungen aus dem Mathematischen Seminar der Universit{\"a}t Hamburg}, volume = {4}, pages = {47--72}, year = {1925} }

@book{Reidemeister1932, author = {K. Reidemeister}, title = {Knotentheorie}, publisher = {Springer}, year = {1932} }

@article{Alexander1923, author = {J. W. Alexander}, title = {A lemma on systems of knotted curves}, journal = {Proceedings of the National Academy of Sciences of the United States of America}, volume = {9}, number = {3}, pages = {93--95}, year = {1923} }

@article{Markov1935, author = {A. A. Markov}, title = {{\"U}ber die freie {\"A}quivalenz geschlossener Z{\"o}pfe}, journal = {Matematicheskii Sbornik}, volume = {43}, number = {1}, pages = {73--78}, year = {1935} }

@book{Birman1974, author = {J. S. Birman}, title = {Braids, Links, and Mapping Class Groups}, publisher = {Princeton University Press}, year = {1974} }

@book{Lickorish1997, author = {W. B. R. Lickorish}, title = {An Introduction to Knot Theory}, publisher = {Springer}, year = {1997} }

@article{Alexander1928, author = {J. W. Alexander}, title = {Topological invariants of knots and links}, journal = {Transactions of the American Mathematical Society}, volume = {30}, pages = {275--306}, year = {1928} }

@incollection{Conway1970,
title = {An enumeration of knots and links, and some of their algebraic properties},
editor = {John Leech},
booktitle = {Computational Problems in Abstract Algebra},
publisher = {Pergamon},
pages = {329-358},
year = {1970},
isbn = {978-0-08-012975-4},
doi = {https://doi.org/10.1016/B978-0-08-012975-4.50034-5},
url = {https://www.sciencedirect.com/science/article/pii/B9780080129754500345},
author = {J.H. Conway}
}

@article{Joyce1982, author = {D. Joyce}, title = {A classifying invariant of knots, the knot quandle}, journal = {Journal of Pure and Applied Algebra}, volume = {23}, number = {1}, pages = {37--65}, year = {1982} }

@article{Jones1985, author = {V. F. R. Jones}, title = {A polynomial invariant for knots via von {N}eumann algebras}, journal = {Bulletin of the American Mathematical Society}, volume = {12}, number = {1}, pages = {103--111}, year = {1985} }

@article{Jones1987, author = {V. F. R. Jones}, title = {Hecke algebra representations of braid groups and link polynomials}, journal = {Annals of Mathematics}, volume = {126}, number = {2}, pages = {335--388}, year = {1987} }

@article{FreydYetterHosteLickorishMillettOcneanu1985, author = {P. Freyd and D. Yetter and J. Hoste and W. B. R. Lickorish and K. Millett and A. Ocneanu}, title = {A new polynomial invariant of knots and links}, journal = {Bulletin of the American Mathematical Society}, volume = {12}, number = {2}, pages = {239--246}, year = {1985} }

@article{PrzytyckiTraczyk1987, author = {J. H. Przytycki and P. Traczyk}, title = {Invariants of links of {C}onway type}, journal = {Kobe Journal of Mathematics}, volume = {4}, pages = {115--139}, year = {1987} }

@article{Kauffman1987, author = {L. H. Kauffman}, title = {State models and the {J}ones polynomial}, journal = {Topology}, volume = {26}, number = {3}, pages = {395--407}, year = {1987} }

@article{BirmanWenzl1989, author = {J. S. Birman and H. Wenzl}, title = {Braids, link polynomials and a new algebra}, journal = {Transactions of the American Mathematical Society}, volume = {313}, pages = {249--273}, year = {1989} }

@article{Murakami1987, author = {J. Murakami}, title = {The {K}auffman polynomial of links and representation theory}, journal = {Osaka Journal of Mathematics}, volume = {24}, pages = {745--758}, year = {1987} }

@book{Kauffman1991KnotsPhysics, author = {L. H. Kauffman}, title = {Knots and Physics}, publisher = {World Scientific}, year = {1991} }

@article{Przytycki1991, author = {J. H. Przytycki}, title = {Skein modules of 3-manifolds}, journal = {Bulletin of the Polish Academy of Sciences}, volume = {39}, pages = {91--100}, year = {1991} }

@article{HosteKidwell1990,
  author  = {J. Hoste and M. Kidwell},
  title   = {Dichromatic link invariants},
  journal = {Transactions of the American Mathematical Society},
  volume  = {321},
  number  = {1},
  pages   = {197--229},
  year    = {1990}
}

@article{Turaev1990SkeinSolidTorus, author = {V. G. Turaev}, title = {The {C}onway and {K}auffman modules of a solid torus}, journal = {Journal of Soviet Mathematics}, volume = {52}, pages = {2799--2805}, year = {1990} }

@article{HaringOldenburgLambropoulou2002Handlebodies, author = {R. H{\"a}ring-Oldenburg and S. Lambropoulou}, title = {Knot theory in handlebodies}, journal = {Journal of Knot Theory and Its Ramifications}, volume = {11}, number = {6}, pages = {921--943}, year = {2002} }

@article{Lambropoulou1999HeckeTypeB, author = {S. Lambropoulou}, title = {Knot theory related to generalized and cyclotomic {H}ecke algebras of type \(B\)}, journal = {Journal of Knot Theory and Its Ramifications}, volume = {8}, number = {5}, pages = {621--658}, year = {1999} }

@article{LambropoulouRourke1997Markov3Manifolds, author = {S. Lambropoulou and C. P. Rourke}, title = {Markov's theorem in 3-manifolds}, journal = {Topology and its Applications}, volume = {78}, pages = {95--122}, year = {1997} }

@incollection{Vassiliev1990,
  author    = {V. A. Vassiliev},
  title     = {Cohomology of Knot Spaces},
  booktitle = {Theory of Singularities and Its Applications},
  editor    = {V. I. Arnold},
  series    = {Advances in Soviet Mathematics},
  volume    = {1},
  publisher = {American Mathematical Society},
  address   = {Providence, RI},
  pages     = {23--69},
  year      = {1990},
  doi       = {10.1090/advsov/001/03}
}

@article{Baez1992, author = {J. C. Baez}, title = {Link invariants of finite type and perturbation theory}, journal = {Letters in Mathematical Physics}, volume = {26}, number = {1}, pages = {43--51}, year = {1992} }

@article{BarNatan1995, author = {D. Bar-Natan}, title = {On the {V}assiliev knot invariants}, journal = {Topology}, volume = {34}, pages = {423--472}, year = {1995} }

@article{Gemein1997SingularBraids, author = {B. Gemein}, title = {Singular braids and {M}arkov's theorem}, journal = {Journal of Knot Theory and Its Ramifications}, volume = {6}, number = {4}, pages = {441--454}, year = {1997} }

@article{FennKeymanRourke1998SingularBraid, author = {R. Fenn and E. Keyman and C. P. Rourke}, title = {The singular braid monoid embeds in a group}, journal = {Journal of Knot Theory and Its Ramifications}, volume = {7}, number = {7}, pages = {881--892}, year = {1998} }

@article{ParisRabenda2004, author = {L. Paris and N. Rabenda}, title = {Singular {H}ecke algebras, {M}arkov traces and link invariants}, journal = {Annales de l'Institut Fourier}, volume = {58}, pages = {2413--2443}, year = {2008} }

@article{Kauffman1999, author = {L. H. Kauffman}, title = {Virtual knot theory}, journal = {European Journal of Combinatorics}, volume = {20}, pages = {663--691}, year = {1999} }

@article{CarterKamadaSaito2002, author = {J. S. Carter and S. Kamada and M. Saito}, title = {Stable equivalence of knots on surfaces and virtual knot cobordisms}, journal = {Journal of Knot Theory and Its Ramifications}, volume = {11}, pages = {311--320}, year = {2002} }

@article{Kuperberg2003, author = {G. Kuperberg}, title = {What is a virtual link?}, journal = {Algebraic \& Geometric Topology}, volume = {3}, pages = {587--591}, year = {2003} }

@article{KauffmanLambropoulou2006VirtualBraids, author = {L. H. Kauffman and S. Lambropoulou}, title = {Virtual braids and the \(L\)-moves}, journal = {Journal of Knot Theory and Its Ramifications}, volume = {15}, number = {6}, pages = {773--811}, year = {2006} }

@article{FennRimanyiRourke1997BraidPermutation, author = {R. Fenn and R. Rim{\'a}nyi and C. Rourke}, title = {The braid-permutation group}, journal = {Topology}, volume = {36}, number = {1}, pages = {123--135}, year = {1997}, doi = {10.1016/0040-9383(95)00072-0} }

@article{Kamada2007WeldedBraids,
  author  = {S. Kamada},
  title   = {Braid Presentation of Virtual Knots and Welded Knots},
  journal = {Osaka Journal of Mathematics},
  volume  = {44},
  number  = {2},
  pages   = {441--458},
  year    = {2007}
}

@article{Damiani2017LoopBraidGroups,
  author  = {C. Damiani},
  title   = {A Journey Through Loop Braid Groups},
  journal = {Expositiones Mathematicae},
  volume  = {35},
  number  = {3},
  pages   = {252--285},
  year    = {2017}
}

@article{Satoh2000,
  author  = {S. Satoh},
  title   = {Virtual Knot Presentation of Ribbon Torus-Knots},
  journal = {Journal of Knot Theory and Its Ramifications},
  volume  = {9},
  number  = {4},
  pages   = {531--542},
  year    = {2000}
}

@article{Manturov2009FreeKnots,
  author  = {V. O. Manturov},
  title   = {On Free Knots and Links},
  journal = {arXiv preprint},
  eprint  = {0902.0127},
  archivePrefix = {arXiv},
  primaryClass = {math.GT},
  year    = {2009}
}

@article{Turaev2012, author = {V. Turaev}, title = {Knotoids}, journal = {Osaka Journal of Mathematics}, volume = {49}, pages = {195--223}, year = {2012} }

@article{GugumcuLambropoulou2021Braidoids, author = {N. G{\"u}g{\"u}mc{\"u} and S. Lambropoulou}, title = {Braidoids}, journal = {Israel Journal of Mathematics}, volume = {242}, pages = {955--995}, year = {2021} }

@article{GoundaroulisEtAl2017BondedKnotoids, author = {D. Goundaroulis and N. G{\"u}g{\"u}mc{\"u} and S. Lambropoulou and J. Dorier and A. Stasiak and L. H. Kauffman}, title = {Topological models for open-knotted protein chains using the concepts of knotoids and bonded knotoids}, journal = {Polymers}, volume = {9}, number = {9}, pages = {444}, year = {2017} }

@article{DorierEtAl2018KnotoidProtein, author = {J. Dorier and D. Goundaroulis and F. Benedetti and A. Stasiak}, title = {Knotoid: a tool to study the entanglement of open protein chains using the concept of knotoids}, journal = {Bioinformatics}, volume = {34}, number = {19}, pages = {3402--3404}, year = {2018} }

@article{Hanaki2010, author = {R. Hanaki}, title = {Pseudo diagrams of links, links and spatial graphs}, journal = {Osaka Journal of Mathematics}, volume = {47}, pages = {863--883}, year = {2010} }

@article{Dye2010, author = {H. A. Dye}, title = {Pseudo knots and an obstruction to cosmetic crossings}, journal = {Journal of Knot Theory and Its Ramifications}, volume = {26}, pages = {1750022}, year = {2017} }

@article{HenrichEtAl2013Pseudoknots, author = {A. Henrich and R. Hoberg and S. Jablan and L. Johnson and E. Minten and L. Radovic}, title = {The theory of pseudoknots}, journal = {Journal of Knot Theory and Its Ramifications}, volume = {22}, number = {7}, pages = {1350032}, year = {2013} }

@article{BardakovJablanWang2016PseudoBraids, author = {V. Bardakov and S. Jablan and H. Wang}, title = {Monoid and group of pseudo braids}, journal = {Journal of Knot Theory and Its Ramifications}, volume = {25}, number = {9}, pages = {1641002}, year = {2016} }

@article{Diamantis2023PseudoSingularSolidTorus, author = {I. Diamantis}, title = {Pseudo Links and Singular Links in the Solid Torus}, journal = {Communications in Mathematics}, volume = {31}, number = {1}, year = {2023} }

@article{Diamantis2021PseudoHandlebodies, author = {I. Diamantis}, title = {Pseudo Links in Handlebodies}, journal = {Bulletin of the Hellenic Mathematical Society}, volume = {65}, pages = {17--34}, year = {2021} }

@misc{Diamantis2026HOMFLYPTPseudoResolution, author = {I. Diamantis}, title = {A HOMFLYPT-type invariant for pseudo links via a resolution in {H}ecke algebras}, eprint = {2605.01026}, archivePrefix = {arXiv}, primaryClass = {math.GT}, note = {Submitted}, year = {2026} }

@misc{AicardiJuyumaya2000BraidsTies, author = {F. Aicardi and J. Juyumaya}, title = {An algebra involving braids and ties}, note = {Preprint ICTP, IC/2000/179, Trieste}, year = {2000} }

@article{AicardiJuyumaya2016TiedLinks, author = {F. Aicardi and J. Juyumaya}, title = {Tied links}, journal = {Journal of Knot Theory and Its Ramifications}, volume = {25}, number = {9}, pages = {}, year = {2016}, note = {DOI:
10.1142/S02182165164100171} }

@article{AicardiJuyumaya2018TiedSingular, author = {F. Aicardi and J. Juyumaya}, title = {Tied links and invariants for singular links}, journal = {Advances in Mathematics}, volume = {381}, pages = {107629}, year = {2021} }

@article{AicardiJuyumaya2018KauffmanTied, author = {F. Aicardi and J. Juyumaya}, title = {Kauffman type invariants for tied links}, journal = {Mathematische Zeitschrift}, volume = {289}, pages = {567--591}, year = {2018} }

@article{Diamantis2021TiedPseudoLinks, author = {I. Diamantis}, title = {Tied pseudo links \& pseudo knotoids}, journal = {Mediterranean Journal of Mathematics}, volume = {18}, pages = {201}, year = {2021} }

@article{Diamantis2021TiedLinksSettings, author = {I. Diamantis}, title = {Tied links in various topological settings}, journal = {Journal of Knot Theory and Its Ramifications}, volume = {30}, number = {7}, pages = {2150046}, year = {2021} }

@article{Kauffman1989Graphs, author = {L. H. Kauffman}, title = {Invariants of graphs in three-space}, journal = {Transactions of the American Mathematical Society}, volume = {311}, number = {2}, pages = {697--710}, year = {1989} }

@incollection{HenrichKauffman2017TangleInsertion, author = {A. K. Henrich and L. H. Kauffman}, title = {Tangle insertion invariants for pseudoknots, singular knots, and rigid vertex spatial graphs}, booktitle = {Knots, Links, Spatial Graphs, and Algebraic Invariants}, series = {Contemporary Mathematics}, volume = {689}, publisher = {American Mathematical Society}, pages = {61--83}, year = {2017} }

@article{DiamantisKauffmanLambropoulou2025Bonded, author = {I. Diamantis and L. H. Kauffman and S. Lambropoulou}, title = {Topology and Algebra of Bonded Knots and Braids}, journal = {Mathematics}, volume = {13}, pages = {3260}, year = {2025}, doi = {10.3390/math13203260} }

@article{Lambropoulou2007LMoves, author = {S. Lambropoulou}, title = {\(L\)-moves and {M}arkov theorems}, journal = {Journal of Knot Theory and Its Ramifications}, volume = {16}, pages = {1459--1468}, year = {2007} }

@article{DiamantisLambropoulouMahmoudi2024ToroidalPseudoKnots, author = {I. Diamantis and S. Lambropoulou and S. Mahmoudi}, title = {From Annular to Toroidal Pseudo Knots}, journal = {Symmetry}, volume = {16}, number = {10}, pages = {1360}, year = {2024} }

@article{DiamantisLambropoulouMahmoudi2026DoublyPeriodicTangles, author = {I. Diamantis and S. Lambropoulou and S. Mahmoudi}, title = {Equivalence of Doubly Periodic Tangles}, journal = {Mathematics}, volume = {14}, number = {6}, pages = {1071}, year = {2026} }

@incollection{DiamantisLambropoulouMahmoudi2024DoublyPeriodicPseudoTangles, author = {I. Diamantis and S. Lambropoulou and S. Mahmoudi}, title = {The Theory of Doubly Periodic Pseudo Tangles}, booktitle = {Contemporary Mathematics}, publisher = {American Mathematical Society}, note = {Accepted for publication. arXiv:2412.16808 [math.GT]}, year = {2024} }

@misc{DiamantisLambropoulouMahmoudi2025ToroidalKnotoids, author = {I. Diamantis and S. Lambropoulou and S. Mahmoudi}, title = {From Annular to Toroidal Knotoids and Their Universal Bracket Polynomials}, eprint = {2509.05014}, archivePrefix = {arXiv}, primaryClass = {math.GT}, year = {2025} }

@article{Bataineh2020,
  author  = {K. Bataineh},
  title   = {Stuck Knots},
  journal = {Symmetry},
  volume  = {12},
  number  = {9},
  pages   = {1558},
  year    = {2020}
}

@article{CenicerosElhamdadiKomissarLahrani2024,
  author  = {J. Ceniceros and M. Elhamdadi and J. Komissar and H. Lahrani},
  title   = {RNA Foldings and Stuck Knots},
  journal = {Communications of the Korean Mathematical Society},
  volume  = {39},
  number  = {1},
  pages   = {223--245},
  year    = {2024}
}

@article{CenicerosElhamdadiMagillRosario2023,
  author  = {J. Ceniceros and M. Elhamdadi and B. Magill and G. Rosario},
  title   = {RNA Foldings, Oriented Stuck Knots and State Sum Invariants},
  journal = {Journal of Mathematical Physics},
  volume  = {64},
  pages   = {031702},
  year    = {2023}
}

@article{Diamantis2026StuckKnots, author = {I. Diamantis}, title = {Stuck Knots: Rigidity, Invariants, and Unsticking Distance}, journal = {Topology and its Applications}, volume = {390}, pages = {109891}, year = {2026}, doi = {10.1016/j.topol.2026.109891} }

@phdthesis{Diamantis2015PhD, author = {I. Diamantis}, title = {The HOMFLYPT Skein Module of the Lens Spaces \(L(p,1)\)}, school = {National Technical University of Athens}, address = {Athens, Greece}, year = {2015}, doi = {10.26240/heal.ntua.1824} }

@article{DiamantisLambropoulou2015BraidEquivalence3Manifolds, author = {I. Diamantis and S. Lambropoulou}, title = {Braid Equivalence in 3-manifolds with Rational Surgery Description}, journal = {Topology and its Applications}, volume = {194}, pages = {269--295}, year = {2015}, doi = {10.1016/j.topol.2015.08.009} }

@article{DiamantisLambropoulou2016NewBasisHOMFLYPTST, author = {I. Diamantis and S. Lambropoulou}, title = {A New Basis for the HOMFLYPT Skein Module of the Solid Torus}, journal = {Journal of Pure and Applied Algebra}, volume = {220}, number = {2}, pages = {577--605}, year = {2016}, doi = {10.1016/j.jpaa.2015.06.014} }

@article{DiamantisLambropoulouPrzytycki2016TopologicalSteps, author = {I. Diamantis and S. Lambropoulou and J. H. Przytycki}, title = {Topological Steps Toward the HOMFLYPT Skein Module of the Lens Spaces \(L(p,1)\) via Braids}, journal = {Journal of Knot Theory and Its Ramifications}, volume = {25}, number = {14}, year = {2016}, doi = {10.1142/S021821651650084X} }

@incollection{DiamantisLambropoulou2017BraidApproach, author = {I. Diamantis and S. Lambropoulou}, title = {The Braid Approach to the HOMFLYPT Skein Module of the Lens Spaces \(L(p,1)\)}, booktitle = {Algebraic Modeling of Topological and Computational Structures and Applications}, series = {Springer Proceedings in Mathematics \& Statistics}, volume = {219}, publisher = {Springer}, year = {2017}, doi = {10.1007/978-3-319-68103-0_7} }

@incollection{Diamantis2019AlternativeBasisST, author = {I. Diamantis}, title = {An Alternative Basis for the Kauffman Bracket Skein Module of the Solid Torus via Braids}, booktitle = {Knots, Low-Dimensional Topology and Applications}, series = {Springer Proceedings in Mathematics \& Statistics}, volume = {284}, publisher = {Springer}, address = {Cham}, year = {2019}, doi = {10.1007/978-3-030-16031-9_16} }

@article{Diamantis2019KBSMGenus2Handlebody, author = {I. Diamantis}, title = {The Kauffman Bracket Skein Module of the Handlebody of Genus 2 via Braids}, journal = {Journal of Knot Theory and Its Ramifications}, volume = {28}, number = {13}, pages = {1940020}, year = {2019}, doi = {10.1142/S0218216519400200} }

@article{DiamantisLambropoulou2019ImportantStep, author = {I. Diamantis and S. Lambropoulou}, title = {An Important Step for the Computation of the HOMFLYPT Skein Module of the Lens Spaces \(L(p,1)\) via Braids}, journal = {Journal of Knot Theory and Its Ramifications}, volume = {28}, number = {11}, pages = {1940007}, year = {2019}, doi = {10.1142/S0218216519400078} }

@article{Diamantis2021HOMFLYPTSubmodulesLens, author = {I. Diamantis}, title = {HOMFLYPT Skein Sub-modules of the Lens Spaces \(L(p,1)\)}, journal = {Topology and its Applications}, volume = {301}, pages = {107500}, year = {2021} }

@article{Diamantis2024KBSMLensUnoriented, author = {I. Diamantis}, title = {The Kauffman Bracket Skein Module of the Lens Spaces via Unoriented Braids}, journal = {Communications in Contemporary Mathematics}, volume = {26}, number = {2}, pages = {2250076}, year = {2024} }

@article{Diamantis2024KBSMS1S2, author = {I. Diamantis}, title = {The Kauffman Bracket Skein Module of \(S^1\times S^2\) via Braids}, journal = {Axioms}, volume = {13}, number = {9}, pages = {617}, year = {2024} }

@article{Diamantis2025HOMFLYPTS1S2, author = {I. Diamantis}, title = {On the HOMFLYPT Skein Module of \(S^1 \times S^2\) via Braids}, journal = {Acta Universitatis Sapientiae, Mathematica}, volume = {17}, pages = {24}, year = {2025}, doi = {10.1007/s44426-025-00025-9} }

@incollection{Diamantis2025SurveySkeinModules, author = {I. Diamantis}, title = {A Survey on Skein Modules via Braids}, booktitle = {Algebraic Structures in Knot Theory}, series = {Contemporary Mathematics}, volume = {827}, publisher = {American Mathematical Society}, address = {Providence, RI}, pages = {49--85}, year = {2025} }

@article{Diamantis2022KnotoidsTorus, author = {I. Diamantis}, title = {Knotoids, Pseudo Knotoids, Braidoids and Pseudo Braidoids on the Torus}, journal = {Communications of the Korean Mathematical Society}, volume = {37}, number = {4}, pages = {1221--1248}, year = {2022} }

\end{document}